\documentstyle[%draft,%showkeys,
amssymb,leqno]{amsart}
\newtheorem{thm}{Theorem}[section]
\newtheorem{lem}[thm]{Lemma}
\newtheorem{prop}[thm]{Proposition}

\theoremstyle{definition}
\newtheorem{defn}[thm]{Definition}

\theoremstyle{remark}
\newtheorem{rem}[thm]{Remark}

\renewcommand{\theequation}{\thesection.\arabic{equation}}
\newcommand{\defref}[1]{Definition~\ref{#1}}

\newcommand{\lemref}[1]{Lemma~\ref{#1}}
\newcommand{\appref}[1]{Appendix~\ref{#1}}
\newcommand{\figref}[1]{Figure~\ref{#1}}

\newcommand{\lam}{\lambda} 
\newcommand{\sig}{\sigma}
\newcommand{\calC}{{\cal C}}
\newcommand{\calH}{{\cal H}}

\newcommand{\I}{{\mathrm I}}
\newcommand{\II}{{\mathrm II}}
\newcommand{\Integer}{{\Bbb Z}}

\newcommand{\Complex}{{\Bbb C}}
\newcommand{\Bethe}{{\frak B}}
\newcommand{\End}{\operatorname{End}}
\newcommand{\Id}{\operatorname{Id}}

\newcommand{\tensor}{\otimes}

\renewcommand{\Im}{\operatorname{Im}}
\newcommand{\s}{\sharp}

\newcommand{\lamcirc}{{\overset{\circ}{\lam}}}
\newcommand{\faceW}[5]{W\left(\left.\begin{matrix}#1&#2\\#3&#4%
\end{matrix}\right|#5\right)}
\begin{document}

\title[Difference equations and Bethe vectors]
{A system of difference equations \\
with elliptic coefficients \\
and Bethe vectors}
\author[Takashi Takebe]{Takashi Takebe\\
        Department of Mathematical Sciences\\
        the University of Tokyo, Komaba 3-8-1, Meguro-ku\\
        Tokyo, 153 Japan}
\address{Department of Mathematical Sciences\\
        the University of Tokyo, Hongo 7-3-1, Bunkyo-ku\\
        Tokyo, 113 Japan}
\curraddr{Department of Mathematics, the University of California,
        Berkeley, CA94720, USA}
\email{takebe@@math.berkeley.edu, takebe@@ms.u-tokyo.ac.jp}
\thanks{Present address: Department of Mathematics, the University of
        California, Berkeley, CA94720, U.S.A. (till August 1997).}

\date{March, 1996}
%
%\keywords{}
\subjclass{Primary
           17B37, % Quantum groups and related deformations,
           Secondary
           39A12, % Discrete version of topics in analysis,
           81R50, % Quantum groups and related algebraic methods,
           82B23  % Exactly solvable models
}
\maketitle
%%%%%%%%%%%%%%%%%%%%%%%%%%%%%%%%%%%%%%%%%%%%%%%%%%%%%%%%%%%%%%%%%%%%%%%%
\centerline{e-mail:
           {\tt takebe@@math.berkeley.edu, takebe@@ms.u-tokyo.ac.jp}}
%%%%%%%%%%%%%%%%%%%%%%%%%%%%%%%%%%%%%%%%%%%%%%%%%%%%%%%%%%%%%%%%%%%%%%%%
\begin{abstract}
An elliptic analogue of the $q$ deformed Knizhnik-Zamolodchikov equations
is introduced. A solution is given in the form of a Jackson-type integral
of Bethe vectors of the XYZ-type spin chains.
\end{abstract}
%%%%%%%%%%%%%%%%%%%%%%%%%%%%%%%%%%%%%%%%%%%%%%%%%%%%%%%%%%%%%%%%%%%%%%%%
%
%\overfullrule=0pt
%%%%%%%%%%%%%%%%%%%%%%%%%%%%%%%%%%%%%%%%%%%%%%%%%%%%%%%%%%%%%%%%%%%%%%%%
%\setlength{\baselineskip}{2\baselineskip} % for double spacing.
%\tableofcontents
\section*{Introduction}

In this paper we introduce a holonomic system of difference equations
associated to elliptic $R$ matrices and give its solution in the form of
Jackson-type integral, following Reshetikhin's idea \cite{resh:92} for the
trigonometric $R$ matrices.

Reshetikhin constructed a solution to the $q$-deformed
Knizhnik-Zamolodchikov equations \cite{fre-resh} by a Jackson-type
integration of Bethe vectors of the XXZ-type spin chain models. Matsuo
\cite{mat} also found the same kind of formulae from a different
viewpoint.

On the other hand, the Bethe Ansatz method for the spin chain models
associated to the elliptic $R$ matrices has been studied since Baxter
\cite{bax}. Hence a natural question is how to find an elliptic version of
Reshetikhin's approach to the $q$-KZ equation. It turns out that the
argument in \cite{resh:92} can be carried out for the elliptic $R$
matrices as well, except for one point. In contrast to the trigonometric
case, an elliptic spin chain model does not have a unique vacuum vector in
its local state space but a series of ``pseudo-vacua'' which depend
non-trivially on a spectral parameter. This dependence breaks down naive
analogy.

We overcome this difficulty by introducing a ``space of Bethe vectors''
and a boundary operator which shifts a spectral parameter.

In the vertex picture the linear space of Bethe vectors depends on
spectral parameters. Therefore we have to use the IRF picture in order to
interpret the system as a holonomic matrix difference system in the sense
of Aomoto \cite{aom:88}. In this context, our system is described in terms
of representation of Felder's elliptic quantum groups \cite{fel:94}.

\medskip
This paper is organized as follows. In the first section we recall several
facts related to the elliptic $R$ matrices and introduce a space of Bethe
vectors and a boundary operator. We define a system of difference
equations in the next section and show its holonomicity. The third section
is devoted to construction of a solution of this system by a Jackson-type
integral of Bethe vectors.

\subsection*{Acknowledgements}
The author expresses his gratitude to Nicolai Reshetikhin for comments and
discussions, and to the Department of Mathematics of the University of
California at Berkeley for hospitality.  He is supported by Postdoctoral
Fellowship for Research abroad of Japan Society for the Promotion of
Science.

\section{Space of Bethe vectors}
\setcounter{equation}{0}

In this section we define a space of Bethe vectors of an elliptic spin
chain and linear operators acting on this space. The state space $\calH$
of a finite XYZ-type spin chain is defined to be a tensor product of the
local state spaces:
\begin{equation}
    \calH := V^{l_1} \tensor V^{l_2} \tensor \cdots \tensor V^{l_N},
\label{def:state-space}
\end{equation}
where $V^{l_j}$ ($j = 1, \ldots, N$) are the spin $l_j$ representation
spaces of the Sklyanin algebra. (See \appref{skl-alg} for a review of the
Sklyanin algebra and its representations.) We will introduce a direct sum
of subspaces of $\calH$ depending on spectral parameters which we will
call a space of Bethe vectors. On this subspace act not only the $R$
matrices but also a boundary operator which shifts the spectral
parameters.

We assume that the parameter $\eta$ which determines the structure
constants of the Sklyanin algebra is a rational number, $\eta = r'/r$.
This assumption is necessary to consider analytic solutions of the system
later.

\subsection*{Space of Bethe vectors}

Baxter \cite{bax} introduced vectors which intertwine vertex-type
Boltzmann weight of the eight-vertex model and IRF-type Boltzmann
weight. They are generalized to higher spin cases. (cf.\ \cite{djkmo},
\cite{take}\footnote{The normalization in previous papers \cite{take} by
the author has been improved. Several complicated factors are absent
now. Vectors $\phi_{\lam,\lam'}(u)$ in those papers correspond to the
following $\phi_{\lam',\lam}(u)$ in the present paper.})

\begin{defn}
\label{def:int-vec}
An {\em $($outgoing$)$ intertwining vector}, 
$\phi^{(l)}_{\lam,\lam'}(u)$, is a vector in the spin $l$ representation
space $V^l = \Theta_{00}^{4l+}$ of the Sklyanin algebra, defined by:
\begin{equation}
\begin{split}
     &\phi^{(l)}_{\lam, \lam'}(u) = \phi^{(l)}_{\lam, \lam'}(u;y) := \\
    :=&
    \prod_{j=1}^{l+m}
    \theta_{10}\left(y + \frac{\lam - u}2 + (2j-l-1)\eta \right)
    \theta_{10}\left(y - \frac{\lam - u}2 - (2j-l-1)\eta \right)
    \times
    \\
    \times
    &\prod_{j=1}^{l-m}
    \theta_{10}\left(y + \frac{\lam'+ u}2 + (2j-l-1)\eta \right)
    \theta_{10}\left(y - \frac{\lam'+ u}2 - (2j-l-1)\eta \right),
\end{split}
\end{equation}
where $\lam$, $\lam'$ are parameters satisfying $\lam - \lam' = 4 m \eta$,
$m \in \{-l, -l+1, \ldots, l\}$, and $u$ is a complex parameter called a
{\em spectral parameter}.  We call
\begin{equation}
    \omega^{(l)}_{\lam}(u) = \phi^{(l)}_{\lam,\lam+4l\eta}(u)
\label{def:loc-pseudo-vac}
\end{equation}
a {\em local pseudo-vacuum}.
\end{defn}

\begin{figure}[ht]
\setlength{\unitlength}{0.00083300in}%
\begingroup\makeatletter\ifx\SetFigFont\undefined
% extract first six characters in \fmtname
\def\x#1#2#3#4#5#6#7\relax{\def\x{#1#2#3#4#5#6}}%
\expandafter\x\fmtname xxxxxx\relax \def\y{splain}%
\ifx\x\y   % LaTeX or SliTeX?
\gdef\SetFigFont#1#2#3{%
  \ifnum #1<17\tiny\else \ifnum #1<20\small\else
  \ifnum #1<24\normalsize\else \ifnum #1<29\large\else
  \ifnum #1<34\Large\else \ifnum #1<41\LARGE\else
     \huge\fi\fi\fi\fi\fi\fi
  \csname #3\endcsname}%
\else
\gdef\SetFigFont#1#2#3{\begingroup
  \count@#1\relax \ifnum 25<\count@\count@25\fi
  \def\x{\endgroup\@setsize\SetFigFont{#2pt}}%
  \expandafter\x
    \csname \romannumeral\the\count@ pt\expandafter\endcsname
    \csname @\romannumeral\the\count@ pt\endcsname
  \csname #3\endcsname}%
\fi
\fi\endgroup
\begin{picture}(1224,1530)(2989,-1279)
\thicklines
\put(6001,-361){\vector( 0,-1){600}}
\multiput(5401,-361)(48.00000,0.00000){26}{\makebox(6.6667,10.0000){\SetFigFont{10}{12}{rm}.}}
\multiput(6001,-361)(0.00000,109.09091){6}{\line( 0, 1){ 54.545}}
\put(6001,-1261){\makebox(0,0)[b]{\smash{$u$}}}
\put(5701,-136){\makebox(0,0)[b]{\smash{$\lambda$}}}
\put(6301,-136){\makebox(0,0)[b]{\smash{$\lambda'$}}}
\end{picture}
\caption{An outgoing intertwining vector $\phi_{\lambda,\lambda'}(u)$.}
\label{fig:int-vec}
\end{figure}
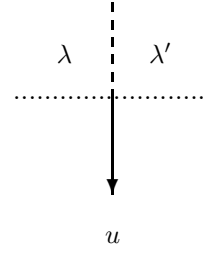

It is easy to see that generically 
$\{ \phi_{\lam+4m\eta,\lam}(u)\}_{m = -l, -l+1, \ldots, l}$
is a basis of $V^l$.  Graphically an intertwining vector is denoted 
as in \figref{fig:int-vec}.

\begin{defn}
A {\em path vector} is a vector in $\calH$ defined as follows:
\begin{multline}
    |a_0, a_1, \ldots, a_N; z_1, \ldots, z_N; \lamcirc \rangle
    :=
    \\
    :=
    \phi^{(l_1)}_{\lam_0, \lam_1}(z_1) \tensor 
    \phi^{(l_2)}_{\lam_1, \lam_2}(z_2) \tensor 
    \cdots \tensor
    \phi^{(l_N)}_{\lam_{N-1}, \lam_N}(z_N) 
\end{multline}
where $a_0, a_1, \ldots, a_N$ are integers, satisfying 
the {\em admissibility condition},
\begin{equation}
    a_{j-1} - a_j \in \{-2l_j, -2l_j+2, \ldots, 2l_j-2, 2l_j\},
\label{admissible}
\end{equation}
$z_j$'s are complex parameters, 
$\lam_n = \lamcirc + 2a_n\eta$ and $\lamcirc$
is a fixed parameter which we omit hereafter unless necessary.
We call a path vector
\begin{equation}
\begin{split}
    \Omega^{l_1, \ldots, l_N}_a(z_1, \ldots, z_N) &= 
    | a_0, a_1, \ldots, a_N ; z_1, \ldots, z_N \rangle \\
    &=
    \omega^{l_1}_{\lam_0}(z_1) \tensor
    \omega^{l_2}_{\lam_1}(z_2)
    \tensor \cdots \tensor
    \omega^{l_N}_{\lam_{N-1}}(z_N)
\end{split}
\label{def:pseudo-vac}
\end{equation}
with $a_0 = a$, $a_j = a + 2(l_1 + \cdots + l_j) \eta$, 
$\lam_j = \lamcirc + 2 a_j \eta$, a {\em $($global$)$ pseudo-vacuum}. 
\end{defn}

We often denote $\Omega^{l_1, \ldots, l_N}_a(z_1, \ldots, z_N)$ 
by $\Omega_a(z_1, \ldots, z_N)$ for simplicity.

Since $\eta = r'/r$ is a rational number,
\begin{equation}
    | a_0 + r, a_1 + r, \ldots, a_N + r ; z_1, \ldots, z_N \rangle
    =
    | a_0    , a_1    , \ldots, a_N     ; z_1, \ldots, z_N \rangle,
\label{periodicity-of-path}
\end{equation}
because of the periodicity of theta functions.

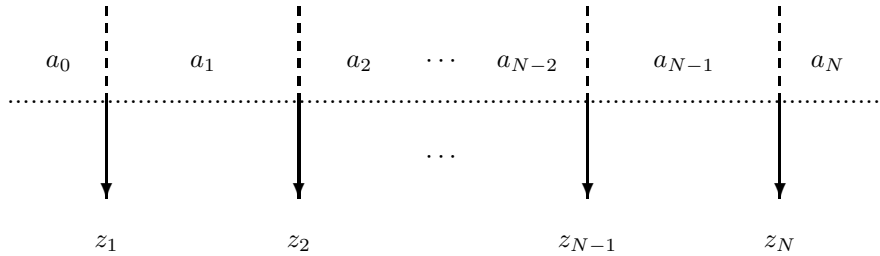
\begin{figure}[ht]
\setlength{\unitlength}{0.00083300in}%
\begingroup\makeatletter\ifx\SetFigFont\undefined
% extract first six characters in \fmtname
\def\x#1#2#3#4#5#6#7\relax{\def\x{#1#2#3#4#5#6}}%
\expandafter\x\fmtname xxxxxx\relax \def\y{splain}%
\ifx\x\y   % LaTeX or SliTeX?
\gdef\SetFigFont#1#2#3{%
  \ifnum #1<17\tiny\else \ifnum #1<20\small\else
  \ifnum #1<24\normalsize\else \ifnum #1<29\large\else
  \ifnum #1<34\Large\else \ifnum #1<41\LARGE\else
     \huge\fi\fi\fi\fi\fi\fi
  \csname #3\endcsname}%
\else
\gdef\SetFigFont#1#2#3{\begingroup
  \count@#1\relax \ifnum 25<\count@\count@25\fi
  \def\x{\endgroup\@setsize\SetFigFont{#2pt}}%
  \expandafter\x
    \csname \romannumeral\the\count@ pt\expandafter\endcsname
    \csname @\romannumeral\the\count@ pt\endcsname
  \csname #3\endcsname}%
\fi
\fi\endgroup
\begin{picture}(5470,1539)(2077,-1288)
\thicklines
\put(3001,-361){\vector( 0,-1){600}}
\multiput(3001,-361)(0.00000,109.09091){6}{\line( 0, 1){ 54.545}}
\put(4201,-361){\vector( 0,-1){600}}
\multiput(4201,-361)(0.00000,109.09091){6}{\line( 0, 1){ 54.545}}
\put(6001,-361){\vector( 0,-1){600}}
\multiput(6001,-361)(0.00000,109.09091){6}{\line( 0, 1){ 54.545}}
\put(7201,-361){\vector( 0,-1){600}}
\multiput(7201,-361)(0.00000,109.09091){6}{\line( 0, 1){ 54.545}}
\multiput(2401,-361)(47.78761,0.00000){114}{\makebox(6.6667,10.0000){\SetFigFont{10}{12}{rm}.}}
\put(2701,-136){\makebox(0,0)[b]{\smash{$a_0$}}}
\put(3601,-136){\makebox(0,0)[b]{\smash{$a_1$}}}
\put(4576,-136){\makebox(0,0)[b]{\smash{$a_2$}}}
\put(7501,-136){\makebox(0,0)[b]{\smash{$a_N$}}}
\put(6601,-136){\makebox(0,0)[b]{\smash{$a_{N-1}$}}}
\put(5626,-136){\makebox(0,0)[b]{\smash{$a_{N-2}$}}}
\put(3001,-1261){\makebox(0,0)[b]{\smash{$z_1$}}}
\put(4201,-1261){\makebox(0,0)[b]{\smash{$z_2$}}}
\put(6001,-1261){\makebox(0,0)[b]{\smash{$z_{N-1}$}}}
\put(7201,-1261){\makebox(0,0)[b]{\smash{$z_N$}}}
\put(5101,-736){\makebox(0,0)[b]{\smash{$\cdots$}}}
\put(5101,-136){\makebox(0,0)[b]{\smash{$\cdots$}}}
\end{picture}
\caption{Path vector 
$|a_0, a_1, \ldots, a_N; z_1, \ldots, z_N \rangle$.}
\label{fig:path-vec}
\end{figure}

A graphical notation for a path vector 
$|a_0, a_1, \ldots, a_N; z_1, \ldots, z_N \rangle$ is shown in
\figref{fig:path-vec}. Though we should write $\lambda_n$ in this figure
instead of $a_n$ for the consistency with \figref{fig:int-vec}, we use
$a_n$ for later convenience.

\begin{defn}
The {\em space of Bethe vectors}
$\Bethe^{l_1, \ldots, l_N}_{z_1, \ldots, z_N}$ with spectral
parameters $\vec z =(z_1, z_2, \ldots, z_N)$ is a space of functions
\begin{equation*}
    f: \Integer / r \Integer \owns \nu \mapsto f(\nu) \in \calH
\end{equation*}
of the form 
\begin{equation*}
    f(\nu) = \sum_{a=0}^{r-1} e^{2\pi i a \nu \eta}
             \sum_{\vec a} \varphi_{\vec a} \,
             | a_0, a_1, \ldots, a_N; z_1, \ldots, z_N \rangle, 
\end{equation*}
where $\varphi_{\vec a}$ are complex numbers and sequences of integers
$\vec a = (a_0, a_1, \ldots, a_N)$ satisfy $a_0 = a_N = a$ and the
admissibility condition \eqref{admissible}. (cf.\
\eqref{periodicity-of-path}.)
\label{def:bethe}
\end{defn}

Generically, 
\begin{equation}
    \left\{\,
    e^{2\pi i a \nu \eta } |a_0, \ldots, a_N; z_1, \ldots, z_N \rangle
    \,\left|
    {{a_0 = a_N = a \in \{0, \ldots, r-1\},}\atop
     {(a_0, \ldots, a_N)\text{ satisfies \eqref{admissible}.}}}
    \right.\right\}
\label{def:bethe-basis}
\end{equation}
is a basis of $\Bethe^{l_1, \ldots, l_N}_{z_1, \ldots, z_N}$. Thus
\begin{equation}
    \dim \Bethe^{l_1, \ldots, l_N}_{z_1, \ldots, z_N}
    =
    r \dim(
    \text{weight zero space of }W^{l_1}\tensor \cdots W^{l_N}),
\end{equation}
where $W^l$ are the spin $l$ (i.e., ($2l+1$)-dimensional) irreducible
representations of the Lie algebra $sl(2, \Complex)$.

\begin{rem}
When $\eta$ is not a rational number (or, exactly speaking, not a
point of finite order on the elliptic curve 
$\Complex / \Integer + \tau \Integer$), the sum should be taken
over all $a \in \Integer$ and $\nu$ is a continuous parameter. The sum
might be considered as a formal series. 
\end{rem}

\subsection*{$R$ matrices}

{\em Elliptic $R$ matrices} $R = R^{l,l'}(u)$ acting on the space $V^l
\tensor V^{l'}$ are constructed by means of the fusion procedure. (cf.\
\cite{cher}, \cite{djkmo}, \cite{hou-zhou}, \cite{take:in-pre}.)  We
recall the following most important properties and refer details to
\cite{take:in-pre}.

\begin{enumerate}
\renewcommand{\labelenumi}{(\roman{enumi})}
\item 
$R^{l,l'}(u)$ is a linear endomorphism of $V^l \tensor V^{l'}$
meromorphically depending on a complex parameter $u$.
(\figref{fig:R-matrix}.)

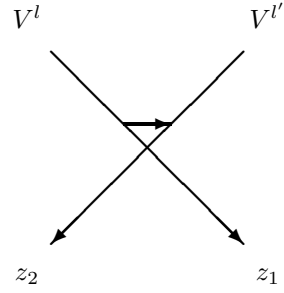
\begin{figure}[t]
\setlength{\unitlength}{0.00083300in}%
\begingroup\makeatletter\ifx\SetFigFont\undefined
% extract first six characters in \fmtname
\def\x#1#2#3#4#5#6#7\relax{\def\x{#1#2#3#4#5#6}}%
\expandafter\x\fmtname xxxxxx\relax \def\y{splain}%
\ifx\x\y   % LaTeX or SliTeX?
\gdef\SetFigFont#1#2#3{%
  \ifnum #1<17\tiny\else \ifnum #1<20\small\else
  \ifnum #1<24\normalsize\else \ifnum #1<29\large\else
  \ifnum #1<34\Large\else \ifnum #1<41\LARGE\else
     \huge\fi\fi\fi\fi\fi\fi
  \csname #3\endcsname}%
\else
\gdef\SetFigFont#1#2#3{\begingroup
  \count@#1\relax \ifnum 25<\count@\count@25\fi
  \def\x{\endgroup\@setsize\SetFigFont{#2pt}}%
  \expandafter\x
    \csname \romannumeral\the\count@ pt\expandafter\endcsname
    \csname @\romannumeral\the\count@ pt\endcsname
  \csname #3\endcsname}%
\fi
\fi\endgroup
\begin{picture}(1815,1755)(3020,-1222)
\thicklines
\put(6601,239){\vector(-1,-1){1200}}
\put(5401,239){\vector( 1,-1){1200}}
\put(5851,-211){\vector( 1, 0){300}}
\put(6751,-1186){\makebox(0,0)[b]{\smash{$z_1$}}}
\put(5251,-1186){\makebox(0,0)[b]{\smash{$z_2$}}}
\put(5251,389){\makebox(0,0)[b]{\smash{$V^l$}}}
\put(6751,389){\makebox(0,0)[b]{\smash{$V^{l'}$}}}
\end{picture}
\caption{$R$ matrix $R^{l,l'}(z_1 - z_2)$.}
\label{fig:R-matrix}
\end{figure}

\item({\em Yang-Baxter equation})
\begin{multline}
    R_{12}^{l ,l' }(z_1 - z_2)
    R_{13}^{l ,l''}(z_1 - z_3)
    R_{23}^{l',l''}(z_2 - z_3) 
    =\\
    =
    R_{23}^{l',l''}(z_2 - z_3) 
    R_{13}^{l ,l''}(z_1 - z_3)
    R_{12}^{l ,l' }(z_1 - z_2)
\label{yb}
\end{multline}
as an endomorphism of $V^l \tensor V^{l'} \tensor V^{l''}$.
(\figref{fig:yb}.)

\begin{figure}[t]
\setlength{\unitlength}{0.00083300in}%
\begingroup\makeatletter\ifx\SetFigFont\undefined
% extract first six characters in \fmtname
\def\x#1#2#3#4#5#6#7\relax{\def\x{#1#2#3#4#5#6}}%
\expandafter\x\fmtname xxxxxx\relax \def\y{splain}%
\ifx\x\y   % LaTeX or SliTeX?
\gdef\SetFigFont#1#2#3{%
  \ifnum #1<17\tiny\else \ifnum #1<20\small\else
  \ifnum #1<24\normalsize\else \ifnum #1<29\large\else
  \ifnum #1<34\Large\else \ifnum #1<41\LARGE\else
     \huge\fi\fi\fi\fi\fi\fi
  \csname #3\endcsname}%
\else
\gdef\SetFigFont#1#2#3{\begingroup
  \count@#1\relax \ifnum 25<\count@\count@25\fi
  \def\x{\endgroup\@setsize\SetFigFont{#2pt}}%
  \expandafter\x
    \csname \romannumeral\the\count@ pt\expandafter\endcsname
    \csname @\romannumeral\the\count@ pt\endcsname
  \csname #3\endcsname}%
\fi
\fi\endgroup
\begin{picture}(5069,2244)(3070,-1822)
\thicklines
\put(5401,239){\vector(-1,-1){1800}}
\put(3601,239){\vector( 1,-1){1800}}
\put(4201,239){\vector( 0,-1){1800}}
\put(4051,-211){\vector( 1, 0){150}}
\put(4426,-586){\makebox(6.6667,10.0000){\SetFigFont{10}{12}{rm}.}}
\put(4426,-586){\makebox(6.6667,10.0000){\SetFigFont{10}{12}{rm}.}}
\put(4351,-511){\vector( 1, 0){300}}
\put(4201,-811){\makebox(6.6667,10.0000){\SetFigFont{10}{12}{rm}.}}
\put(4201,-811){\vector( 1, 0){150}}
\put(8401,239){\vector(-1,-1){1800}}
\put(6601,239){\vector( 1,-1){1800}}
\put(7426,-586){\makebox(6.6667,10.0000){\SetFigFont{10}{12}{rm}.}}
\put(7426,-586){\makebox(6.6667,10.0000){\SetFigFont{10}{12}{rm}.}}
\put(7351,-511){\vector( 1, 0){300}}
\put(7801,-211){\vector( 1, 0){150}}
\put(7651,-811){\vector( 1, 0){150}}
\put(7801,239){\vector( 0,-1){1800}}
\put(3601,314){\makebox(0,0)[b]{\smash{$V^l$}}}
\put(4201,314){\makebox(0,0)[b]{\smash{$V^{l'}$}}}
\put(5401,314){\makebox(0,0)[b]{\smash{$V^{l''}$}}}
\put(6601,314){\makebox(0,0)[b]{\smash{$V^l$}}}
\put(8401,314){\makebox(0,0)[b]{\smash{$V^{l''}$}}}
\put(8401,-1786){\makebox(0,0)[b]{\smash{$z_1$}}}
\put(7801,314){\makebox(0,0)[b]{\smash{$V^{l'}$}}}
\put(6001,-736){\makebox(0,0)[b]{\smash{$=$}}}
\put(7801,-1786){\makebox(0,0)[b]{\smash{$z_2$}}}
\put(6601,-1786){\makebox(0,0)[b]{\smash{$z_3$}}}
\put(5401,-1786){\makebox(0,0)[b]{\smash{$z_1$}}}
\put(4201,-1786){\makebox(0,0)[b]{\smash{$z_2$}}}
\put(3601,-1786){\makebox(0,0)[b]{\smash{$z_3$}}}
\end{picture}
\caption{Yang-Baxter equation.}
\label{fig:yb}
\end{figure}
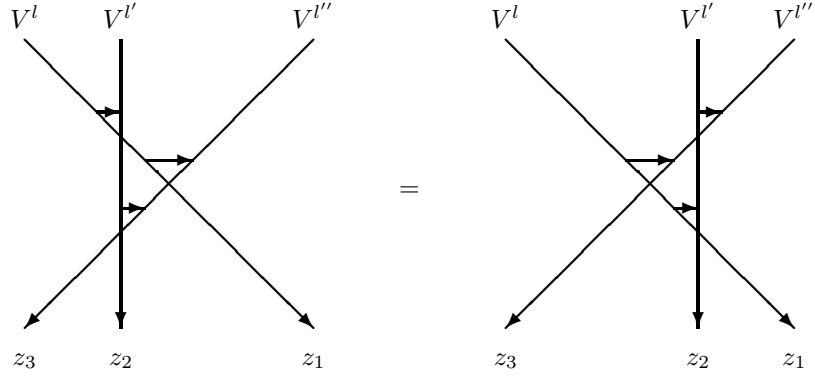

\item({\em Unitarity})
\begin{equation}
    R_{12}^{l,l'}(u - v) R_{21}^{l',l}(v - u) =
    \Id_{V^l\tensor V^{l'}},
\label{unitarity}
\end{equation}
as an endomorphism of $V^l \tensor V^{l'}$.
(\figref{fig:unitarity}.)

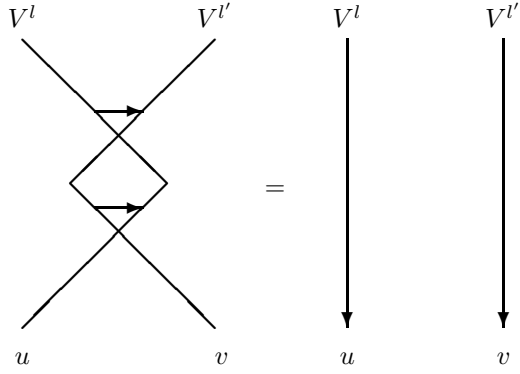
\begin{figure}[t]
\setlength{\unitlength}{0.00083300in}%
\begingroup\makeatletter\ifx\SetFigFont\undefined
% extract first six characters in \fmtname
\def\x#1#2#3#4#5#6#7\relax{\def\x{#1#2#3#4#5#6}}%
\expandafter\x\fmtname xxxxxx\relax \def\y{splain}%
\ifx\x\y   % LaTeX or SliTeX?
\gdef\SetFigFont#1#2#3{%
  \ifnum #1<17\tiny\else \ifnum #1<20\small\else
  \ifnum #1<24\normalsize\else \ifnum #1<29\large\else
  \ifnum #1<34\Large\else \ifnum #1<41\LARGE\else
     \huge\fi\fi\fi\fi\fi\fi
  \csname #3\endcsname}%
\else
\gdef\SetFigFont#1#2#3{\begingroup
  \count@#1\relax \ifnum 25<\count@\count@25\fi
  \def\x{\endgroup\@setsize\SetFigFont{#2pt}}%
  \expandafter\x
    \csname \romannumeral\the\count@ pt\expandafter\endcsname
    \csname @\romannumeral\the\count@ pt\endcsname
  \csname #3\endcsname}%
\fi
\fi\endgroup
\begin{picture}(3184,2226)(3924,-1804)
\thicklines
\put(8401,239){\vector( 0,-1){1800}}
\put(7426,239){\vector( 0,-1){1800}}
\put(6601,239){\line(-1,-1){900}}
\put(6301,-661){\line(-1,-1){900}}
\put(5401,239){\line( 1,-1){900}}
\put(5701,-661){\line( 1,-1){900}}
\put(5851,-811){\vector( 1, 0){300}}
\put(5851,-211){\vector( 1, 0){300}}
\put(5401,314){\makebox(0,0)[b]{\smash{$V^l$}}}
\put(6601,314){\makebox(0,0)[b]{\smash{$V^{l'}$}}}
\put(8401,314){\makebox(0,0)[b]{\smash{$V^{l'}$}}}
\put(7426,314){\makebox(0,0)[b]{\smash{$V^l$}}}
\put(7426,-1786){\makebox(0,0)[b]{\smash{$u$}}}
\put(8401,-1786){\makebox(0,0)[b]{\smash{$v$}}}
\put(5401,-1786){\makebox(0,0)[b]{\smash{$u$}}}
\put(6601,-1786){\makebox(0,0)[lb]{\smash{$v$}}}
\put(6976,-736){\makebox(0,0)[b]{\smash{$=$}}}
\end{picture}
\caption{Unitarity.}
\label{fig:unitarity}
\end{figure}

\item
$R^{l,l'}(u)$ acts on the intertwining vectors as follows:
\begin{multline}
    R^{l,l'}(u-v) 
    \phi_{\lam,\lam'}^{(l)}(u) \tensor \phi_{\lam',\mu}^{(l')}(v)
    =\\
    =
    \sum_{\mu'}
    \faceW{\lam}{\lam'}{\mu'}{\mu}{u-v}
    \phi_{\mu',\mu}^{(l)}(u) \tensor \phi_{\lam,\mu'}^{(l')}(v),
\label{R-on-int-vec}
\end{multline}
where the sum on the right hand side is taken over $\mu'$ satisfying
$\mu' - \mu = 4 m \eta$ ($m \in \{-l, -l+1, \ldots, l\}$) and 
$\lam - \mu' = 4 m' \eta$ ($m \in \{-l', -l'+1, \ldots, l'\}$).
Scalar factors $W$ are the Boltzmann weights of the IRF-type model.
(\figref{fig:R-on-int}. In the figure $W$ is denoted by a crossing of
dashed lines.)

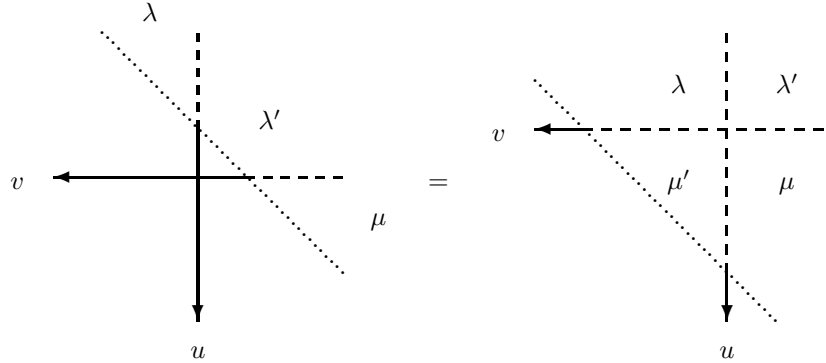
\begin{figure}[t]
\setlength{\unitlength}{0.00083300in}%
\begingroup\makeatletter\ifx\SetFigFont\undefined
% extract first six characters in \fmtname
\def\x#1#2#3#4#5#6#7\relax{\def\x{#1#2#3#4#5#6}}%
\expandafter\x\fmtname xxxxxx\relax \def\y{splain}%
\ifx\x\y   % LaTeX or SliTeX?
\gdef\SetFigFont#1#2#3{%
  \ifnum #1<17\tiny\else \ifnum #1<20\small\else
  \ifnum #1<24\normalsize\else \ifnum #1<29\large\else
  \ifnum #1<34\Large\else \ifnum #1<41\LARGE\else
     \huge\fi\fi\fi\fi\fi\fi
  \csname #3\endcsname}%
\else
\gdef\SetFigFont#1#2#3{\begingroup
  \count@#1\relax \ifnum 25<\count@\count@25\fi
  \def\x{\endgroup\@setsize\SetFigFont{#2pt}}%
  \expandafter\x
    \csname \romannumeral\the\count@ pt\expandafter\endcsname
    \csname @\romannumeral\the\count@ pt\endcsname
  \csname #3\endcsname}%
\fi
\fi\endgroup
\begin{picture}(5154,2226)(4029,-1804)
\thicklines
\put(6001,-661){\vector(-1, 0){1200}}
\multiput(6001,-661)(109.09091,0.00000){6}{\line( 1, 0){ 54.545}}
\multiput(5701,-361)(0.00000,109.09091){6}{\line( 0, 1){ 54.545}}
\multiput(5101,239)(34.09091,-34.09091){45}{\makebox(6.6667,10.0000){\SetFigFont{10}{12}{rm}.}}
\put(5701,-361){\vector( 0,-1){1200}}
\multiput(9001,-1261)(0.00000,120.00000){13}{\line( 0, 1){ 60.000}}
\put(9001,-1261){\vector( 0,-1){300}}
\put(8101,-361){\vector(-1, 0){300}}
\multiput(8101,-361)(120.00000,0.00000){13}{\line( 1, 0){ 60.000}}
\multiput(7801,-61)(34.09091,-34.09091){45}{\makebox(6.6667,10.0000){\SetFigFont{10}{12}{rm}.}}
\put(5701,-1786){\makebox(0,0)[b]{\smash{$u$}}}
\put(5401,314){\makebox(0,0)[b]{\smash{$\lambda$}}}
\put(4576,-736){\makebox(0,0)[b]{\smash{$v$}}}
\put(6826,-961){\makebox(0,0)[b]{\smash{$\mu$}}}
\put(9001,-1786){\makebox(0,0)[b]{\smash{$u$}}}
\put(7576,-436){\makebox(0,0)[b]{\smash{$v$}}}
\put(9376,-136){\makebox(0,0)[b]{\smash{$\lambda'$}}}
\put(8701,-736){\makebox(0,0)[b]{\smash{$\mu'$}}}
\put(9376,-736){\makebox(0,0)[b]{\smash{$\mu$}}}
\put(8701,-136){\makebox(0,0)[b]{\smash{$\lambda$}}}
\put(6151,-361){\makebox(0,0)[b]{\smash{$\lambda'$}}}
\put(7201,-736){\makebox(0,0)[b]{\smash{$=$}}}
\end{picture}
\caption{$R$ and IRF weight are intertwined.}
\label{fig:R-on-int}
\end{figure}

\end{enumerate}

Thanks to (iv), the $R$ matrix defines a map between the spaces of Bethe
vectors:
\begin{multline}
    \check R_{j,j+1}(z_j - z_{j+1}) := \\
    :=
    P_{j,j+1} R_{j,j+1}(z_j - z_{j+1}):
    \Bethe^{l_1, \ldots, l_j, l_{j+1}, \ldots, l_N}_{%
            z_1, \ldots, z_j, z_{j+1}, \ldots, z_N}
     \to
    \Bethe^{l_1, \ldots, l_{j+1}, l_j, \ldots, l_N}_{%
            z_1, \ldots, z_{j+1}, z_j, \ldots, z_N}
\label{def:R-check}
\end{multline}
where $P_{j,j+1}$ is a permutation operator of the $j$-th and the
$(j+1)$-st component of the tensor product 
$V^{l_1}\tensor \cdots \tensor V^{l_N}$. With respect to the basis
\eqref{def:bethe-basis}, the $R$ matrix is described by the IRF-type
Boltzmann weight, and thus by a representation of Felder's elliptic
quantum group \cite{fel:94}. For our purpose, the following
special cases of \eqref{R-on-int-vec} are important.

When both $\phi_{\lam,\lam'}^{(l)}(v)$ and $\phi_{\lam',\mu}^{(l')}(u)$
are local pseudo-vacua \eqref{def:loc-pseudo-vac}, i.e., when
$\lam'= \lam+ 4l\eta$ and $\mu = \lam' + 4l'\eta$,
equation \eqref{R-on-int-vec} is simply
\begin{equation}
    R^{l,l'}(u-v) 
    \omega_{\lam        }^{(l )}(u) \tensor 
    \omega_{\lam-4l\eta }^{(l')}(v)
    =
    \omega_{\lam+4l'\eta}^{(l )}(u) \tensor
    \omega_{\lam        }^{(l')}(v).
\label{R-on-vac}
\end{equation}

When $l = 1/2$, the $R$ matrix $R^{1/2,l}$ is expressed as the
$L$ operator \eqref{def:L} through the identification \eqref{rep:identify}:
\begin{equation}
    R^{1/2,l} (u) 
    = \frac{\theta_{11}(2\eta)}{\theta_{11}(u+(2l+1)\eta)}
      \rho^l(L(u+\eta))
    = \sum_{a=0}^3 
      \frac{\theta_{11}(2\eta) W_a^L(u+\eta)}{\theta_{11}(u+(2l+1)\eta)} 
      \sig^a \tensor \rho^{(l)}(S^a).
\label{R=L}
\end{equation}
In particular, $R^{1/2, 1/2}(u)$ is proportional to Baxter's $R$ matrix
\eqref{def:R}. Likewise the intertwining vectors 
$\phi^{(1/2)}_{\lam, \lam\pm 2\eta}(u)$ correspond to the following
vectors in $\Complex^2$
\begin{equation}
    \phi^{(1/2)}_{\lam\pm 2\eta, \lam}(u-\eta)
    = C 
    \begin{pmatrix}
             - \theta_{01}\left((\lam \pm u)/2 ; \tau/2 \right) \\
    \phantom{-}\theta_{00}\left((\lam \pm u)/2 ; \tau/2 \right)
    \end{pmatrix},
\label{int-vec:ident}
\end{equation}
where $C$ is a constant:
\begin{equation*}
    C = e^{- \pi i \tau/8}
    \frac{\theta_{00}(0;\tau)^2 \theta_{01}(0;\tau)\theta_{10}(0;\tau)} 
         {2 \theta_{10}(0;\tau/2) \theta_{01}(0;2\tau)
            \theta_{10}((1+\tau)/4;\tau) \theta_{10}((1-\tau)/4;\tau)}.
\end{equation*}

\begin{figure}[ht]
\setlength{\unitlength}{0.00083300in}%
\begingroup\makeatletter\ifx\SetFigFont\undefined
% extract first six characters in \fmtname
\def\x#1#2#3#4#5#6#7\relax{\def\x{#1#2#3#4#5#6}}%
\expandafter\x\fmtname xxxxxx\relax \def\y{splain}%
\ifx\x\y   % LaTeX or SliTeX?
\gdef\SetFigFont#1#2#3{%
  \ifnum #1<17\tiny\else \ifnum #1<20\small\else
  \ifnum #1<24\normalsize\else \ifnum #1<29\large\else
  \ifnum #1<34\Large\else \ifnum #1<41\LARGE\else
     \huge\fi\fi\fi\fi\fi\fi
  \csname #3\endcsname}%
\else
\gdef\SetFigFont#1#2#3{\begingroup
  \count@#1\relax \ifnum 25<\count@\count@25\fi
  \def\x{\endgroup\@setsize\SetFigFont{#2pt}}%
  \expandafter\x
    \csname \romannumeral\the\count@ pt\expandafter\endcsname
    \csname @\romannumeral\the\count@ pt\endcsname
  \csname #3\endcsname}%
\fi
\fi\endgroup
\begin{picture}(4709,2424)(4104,-1873)
\thicklines
\put(6301,-661){\vector(-1, 0){1500}}
\multiput(5101,-61)(48.00000,0.00000){26}{\makebox(6.6667,10.0000){\SetFigFont{10}{12}{rm}.}}
\multiput(6301,-61)(0.00000,-48.00000){26}{\makebox(6.6667,10.0000){\SetFigFont{10}{12}{rm}.}}
\multiput(6301,-1261)(-48.00000,0.00000){26}{\makebox(6.6667,10.0000){\SetFigFont{10}{12}{rm}.}}
\multiput(6301,-661)(109.09091,0.00000){6}{\line( 1, 0){ 54.545}}
\multiput(5701,-1261)(0.00000,-109.09091){6}{\line( 0,-1){ 54.545}}
\multiput(5701,539)(0.00000,-109.09091){6}{\line( 0,-1){ 54.545}}
\put(5701,-61){\vector( 0,-1){1200}}
\put(8101,-661){\vector(-1, 0){300}}
\multiput(8101,-61)(0.00000,-48.00000){26}{\makebox(6.6667,10.0000){\SetFigFont{10}{12}{rm}.}}
\multiput(8701,-61)(0.00000,-114.28571){11}{\line( 0,-1){ 57.143}}
\multiput(8101,-661)(114.28571,0.00000){11}{\line( 1, 0){ 57.143}}
\put(5401,239){\makebox(0,0)[b]{\smash{$\mu$}}}
\put(6601,239){\makebox(0,0)[b]{\smash{$\lambda'$}}}
\put(6601,-1561){\makebox(0,0)[b]{\smash{$\lambda$}}}
\put(5401,-1561){\makebox(0,0)[b]{\smash{$\mu'$}}}
\put(4651,-736){\makebox(0,0)[b]{\smash{$v$}}}
\put(5626,-961){\makebox(0,0)[rb]{\smash{$u-\eta$}}}
\put(7651,-736){\makebox(0,0)[b]{\smash{$v$}}}
\put(8701,-1411){\makebox(0,0)[b]{\smash{$u-\eta$}}}
\put(7276,-736){\makebox(0,0)[b]{\smash{$=$}}}
\put(9001,-361){\makebox(0,0)[b]{\smash{$\lambda'$}}}
\put(8401,-361){\makebox(0,0)[b]{\smash{$\mu$}}}
\put(8401,-961){\makebox(0,0)[b]{\smash{$\mu'$}}}
\put(9001,-961){\makebox(0,0)[b]{\smash{$\lambda$}}}
\end{picture}
\caption{An element of a twisted $L$ operator
         acts on $\phi_{\lam',\lam}(v)$.}
\label{fig:twisted-L-on-int-vec}
\end{figure}
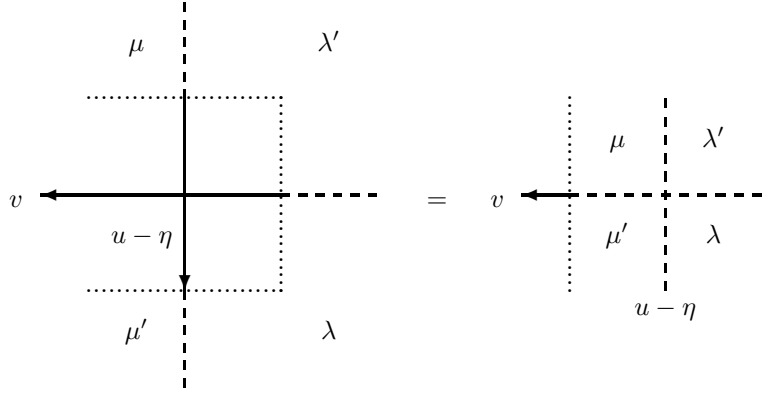

Under this identification, the action of $R^{1/2,l}$ \eqref{R-on-int-vec}
can be stated in the following form which will be used later. Let us
define a matrix of the gauge transformation by
\begin{equation}
    M_{\lam}(u) :=
    \begin{pmatrix}
             - \theta_{01}\left((\lam - u)/2 ; \tau/2\right) &
             - \theta_{01}\left((\lam + u)/2 ; \tau/2\right) \\
    \phantom{-}\theta_{00}\left((\lam - u)/2 ; \tau/2\right) &
    \phantom{-}\theta_{00}\left((\lam + u)/2 ; \tau/2\right)
    \end{pmatrix}
    \begin{pmatrix}
    1 & 0 \\  0 & \theta_{11}(\lam;\tau)^{-1}
    \end{pmatrix},
\label{def:gauge-trans}
\end{equation}
and a twisted $L$ operator by
\begin{equation}
    L_{\lam,\lam'}(u;v) =
    \begin{pmatrix}
    \alpha_{\lam,\lam'}(u;v) &  \beta_{\lam,\lam'}(u;v) \\ 
    \gamma_{\lam,\lam'}(u;v) & \delta_{\lam,\lam'}(u;v)
    \end{pmatrix}
    := M_\lam(u)^{-1} L(u-v) M_{\lam'}(u).
\label{def:twisted-L}
\end{equation}
Then the matrix elements of $L_{\lam,\lam'}(u;v)$ act on an intertwining
vector as follows: 
Let $\lam - \lam' = 4m\eta$ ($m \in \{-l, -l+1, \ldots, l\}$). Then,
\begin{align}
    \alpha_{\lam,\lam'}(u;v) \phi^{(l)}_{\lam',\lam}(v)
    &=
    \frac{\theta_{11}(u - v + 2m\eta) \theta_{11}(\lam + 2(l-m)\eta)}
         {\theta_{11}(\lam) \theta_{11}(2\eta)}
    \phi^{(l)}_{\lam'-2\eta,\lam-2\eta}(v),
\label{alpha-on-int-vec}
\\
    \beta_{\lam,\lam'}(u;v) \phi^{(l)}_{\lam',\lam}(v)
    &=
    \frac{\theta_{11}(u - v + \lam - 2m\eta) \theta_{11}(2(l+m)\eta)}
         {\theta_{11}(\lam) \theta_{11}(\lam-4m\eta) \theta_{11}(2\eta)}
    \phi^{(l)}_{\lam'+2\eta,\lam-2\eta}(v),
\label{beta-on-int-vec}
\\
    \gamma_{\lam,\lam'}(u;v) \phi^{(l)}_{\lam',\lam}(v)
    &=
    \frac{\theta_{11}(u - v - \lam + 2m\eta) \theta_{11}(2(-l+m)\eta)}
         {\theta_{11}(2\eta)}
    \phi^{(l)}_{\lam'-2\eta,\lam+2\eta}(v),
\label{gamma-on-int-vec}
\\
    \delta_{\lam,\lam'}(u;v) \phi^{(l)}_{\lam',\lam}(v)
    &=
    \frac{\theta_{11}(u - v - 2m\eta) \theta_{11}(\lam - 2(l+m)\eta)}
         {\theta_{11}(\lam-4m\eta) \theta_{11}(2\eta)}
    \phi^{(l)}_{\lam'+2\eta,\lam+2\eta}(v).
\label{delta-on-int-vec}
\end{align}
This formulation was found in the context of the quantum inverse
scattering method and the algebraic Bethe Ansatz \cite{takh-fad}. See
also \cite{take}. Note that the column vectors of $M_\lam(u)$ are
essentially outgoing intertwining vectors. Therefore, defining 
the {\em incoming intertwining vectors} $\bar\phi_{\lam,\lam'}(u)$
as row vectors of $M_\lam(u)^{-1}$:
\begin{equation}
    M_\lam(u)^{-1} =
    \begin{pmatrix}
     \bar\phi_{\lam,\lam-2\eta}(u - \eta) \\
     \bar\phi_{\lam,\lam+2\eta}(u - \eta)
    \end{pmatrix},
\label{def:incoming-int-vec}
\end{equation}
and denoting them as in \figref{fig:incoming-int-vec}, 
we can rewrite formulae \eqref{alpha-on-int-vec}--\eqref{delta-on-int-vec}
as in \figref{fig:twisted-L-on-int-vec}. (Exactly speaking, the
normalization here is different from that in \figref{fig:R-on-int}, which
is not essential. In general, incoming intertwining vectors are defined 
as a dual basis of 
$\{ \phi_{\lam+4m\eta,\lam}(u)\}_{m = -l, -l+1, \ldots, l}$. 
See \cite{has:93}.) 

\begin{figure}[ht]
\setlength{\unitlength}{0.00083300in}%
\begingroup\makeatletter\ifx\SetFigFont\undefined
% extract first six characters in \fmtname
\def\x#1#2#3#4#5#6#7\relax{\def\x{#1#2#3#4#5#6}}%
\expandafter\x\fmtname xxxxxx\relax \def\y{splain}%
\ifx\x\y   % LaTeX or SliTeX?
\gdef\SetFigFont#1#2#3{%
  \ifnum #1<17\tiny\else \ifnum #1<20\small\else
  \ifnum #1<24\normalsize\else \ifnum #1<29\large\else
  \ifnum #1<34\Large\else \ifnum #1<41\LARGE\else
     \huge\fi\fi\fi\fi\fi\fi
  \csname #3\endcsname}%
\else
\gdef\SetFigFont#1#2#3{\begingroup
  \count@#1\relax \ifnum 25<\count@\count@25\fi
  \def\x{\endgroup\@setsize\SetFigFont{#2pt}}%
  \expandafter\x
    \csname \romannumeral\the\count@ pt\expandafter\endcsname
    \csname @\romannumeral\the\count@ pt\endcsname
  \csname #3\endcsname}%
\fi
\fi\endgroup
\begin{picture}(1224,1434)(2989,-1273)
\thicklines
\multiput(5401,-661)(48.00000,0.00000){26}{\makebox(6.6667,10.0000){\SetFigFont{10}{12}{rm}.}}
\multiput(6001,-661)(0.00000,-109.09091){6}{\line( 0,-1){ 54.545}}
\put(6001,-61){\vector( 0,-1){600}}
\put(5701,-886){\makebox(0,0)[b]{\smash{$\lambda$}}}
\put(6301,-886){\makebox(0,0)[b]{\smash{$\lambda'$}}}
\put(6001, 89){\makebox(0,0)[b]{\smash{$u$}}}
\end{picture}
\caption{An incoming intertwining vector $\bar\phi_{\lam',\mu}(u)$.}
\label{fig:incoming-int-vec}
\end{figure}
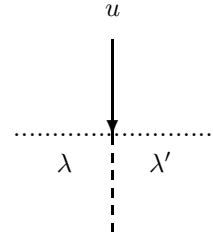

Especially $\alpha$, $\gamma$ and $\delta$ act on a
local pseudo vacuum as follows:
\begin{align}
    \alpha_{\lam,\lam'}(u;v) \omega^{(l)}_{\lam'}(v)
    &= \alpha^l(u - v) \omega^{(l)}_{\lam'-2\eta}(v),
\label{alpha-on-vac}
\\
    \gamma_{\lam,\lam'}(u;v) \omega^{(l)}_{\lam'}(v)
    &= 0
\label{gamma-on-vac}
\\
    \delta_{\lam,\lam'}(u;v) \omega^{(l)}_{\lam'}(v)
    &= \delta^l(u - v) \omega^{(l)}_{\lam'+2\eta}(v),
\label{delta-on-vac}
\end{align}
where
\begin{equation}
\alpha^l(u) = \frac{\theta_{11}(u + 2l\eta)}{\theta_{11}(2\eta)}, \qquad
\delta^l(u) = \frac{\theta_{11}(u - 2l\eta)}{\theta_{11}(2\eta)}.
\label{def:alpha-delta}
\end{equation}

\subsection*{Boundary operator}

We introduce an operator $Z$ which shifts a spectral parameter of a
Bethe vector. Let $\kappa$ and $c$ be fixed complex parameters.

\begin{figure}[ht]
\setlength{\unitlength}{0.00083300in}%
\begingroup\makeatletter\ifx\SetFigFont\undefined
% extract first six characters in \fmtname
\def\x#1#2#3#4#5#6#7\relax{\def\x{#1#2#3#4#5#6}}%
\expandafter\x\fmtname xxxxxx\relax \def\y{splain}%
\ifx\x\y   % LaTeX or SliTeX?
\gdef\SetFigFont#1#2#3{%
  \ifnum #1<17\tiny\else \ifnum #1<20\small\else
  \ifnum #1<24\normalsize\else \ifnum #1<29\large\else
  \ifnum #1<34\Large\else \ifnum #1<41\LARGE\else
     \huge\fi\fi\fi\fi\fi\fi
  \csname #3\endcsname}%
\else
\gdef\SetFigFont#1#2#3{\begingroup
  \count@#1\relax \ifnum 25<\count@\count@25\fi
  \def\x{\endgroup\@setsize\SetFigFont{#2pt}}%
  \expandafter\x
    \csname \romannumeral\the\count@ pt\expandafter\endcsname
    \csname @\romannumeral\the\count@ pt\endcsname
  \csname #3\endcsname}%
\fi
\fi\endgroup
\begin{picture}(5707,3648)(2047,-3397)
\thicklines
\put(4201,-2461){\vector( 0,-1){600}}
\multiput(4201,-2461)(0.00000,109.09091){6}{\line( 0, 1){ 54.545}}
\put(5401,-2461){\vector( 0,-1){600}}
\multiput(5401,-2461)(0.00000,109.09091){6}{\line( 0, 1){ 54.545}}
\put(7201,-2461){\vector( 0,-1){600}}
\multiput(7201,-2461)(0.00000,109.09091){6}{\line( 0, 1){ 54.545}}
\put(3001,-361){\vector( 0,-1){600}}
\multiput(3001,-361)(0.00000,109.09091){6}{\line( 0, 1){ 54.545}}
\put(4201,-361){\vector( 0,-1){600}}
\multiput(4201,-361)(0.00000,109.09091){6}{\line( 0, 1){ 54.545}}
\put(6001,-361){\vector( 0,-1){600}}
\multiput(6001,-361)(0.00000,109.09091){6}{\line( 0, 1){ 54.545}}
\put(7201,-361){\vector( 0,-1){600}}
\multiput(7201,-361)(0.00000,109.09091){6}{\line( 0, 1){ 54.545}}
\put(3001,-2461){\vector( 0,-1){600}}
\multiput(3001,-2461)(0.00000,109.09091){6}{\line( 0, 1){ 54.545}}
\put(4951,-1261){\line( 0,-1){300}}
\put(5251,-1261){\line( 0,-1){300}}
\put(4801,-1411){\line( 1,-1){300}}
\put(5101,-1711){\line( 1, 1){300}}
\multiput(2401,-361)(47.78761,0.00000){114}{\makebox(6.6667,10.0000){\SetFigFont{10}{12}{rm}.}}
\multiput(2401,-2461)(47.78761,0.00000){114}{\makebox(6.6667,10.0000){\SetFigFont{10}{12}{rm}.}}
\put(4801,-2236){\makebox(0,0)[b]{\smash{$a_1$}}}
\put(5776,-2236){\makebox(0,0)[b]{\smash{$a_2$}}}
\put(6826,-2236){\makebox(0,0)[b]{\smash{$a_{N-2}$}}}
\put(4201,-3361){\makebox(0,0)[b]{\smash{$z_1$}}}
\put(5401,-3361){\makebox(0,0)[b]{\smash{$z_2$}}}
\put(7201,-3361){\makebox(0,0)[b]{\smash{$z_{N-1}$}}}
\put(2701,-136){\makebox(0,0)[b]{\smash{$a_0$}}}
\put(3601,-136){\makebox(0,0)[b]{\smash{$a_1$}}}
\put(4576,-136){\makebox(0,0)[b]{\smash{$a_2$}}}
\put(7501,-136){\makebox(0,0)[b]{\smash{$a_N$}}}
\put(6601,-136){\makebox(0,0)[b]{\smash{$a_{N-1}$}}}
\put(5626,-136){\makebox(0,0)[b]{\smash{$a_{N-2}$}}}
\put(3001,-1261){\makebox(0,0)[b]{\smash{$z_1$}}}
\put(4201,-1261){\makebox(0,0)[b]{\smash{$z_2$}}}
\put(6001,-1261){\makebox(0,0)[b]{\smash{$z_{N-1}$}}}
\put(7201,-1261){\makebox(0,0)[b]{\smash{$z_N$}}}
\put(3601,-2236){\makebox(0,0)[b]{\smash{$a_0$}}}
\put(7501,-2236){\makebox(0,0)[b]{\smash{$a_{N-1}$}}}
\put(5101,-136){\makebox(0,0)[b]{\smash{$\cdots$}}}
\put(5101,-811){\makebox(0,0)[b]{\smash{$\cdots$}}}
\put(6301,-2236){\makebox(0,0)[b]{\smash{$\cdots$}}}
\put(6301,-2911){\makebox(0,0)[b]{\smash{$\cdots$}}}
\put(3001,-3361){\makebox(0,0)[b]{\smash{$z_N+\kappa$}}}
\put(2701,-2236){\makebox(0,0)[b]{\smash{$a_{N-1}$}}}
\put(5101,-1561){\makebox(0,0)[b]{\smash{$Z$}}}
\end{picture}
\caption{A boundary operator.}
\label{fig:Z}
\end{figure}
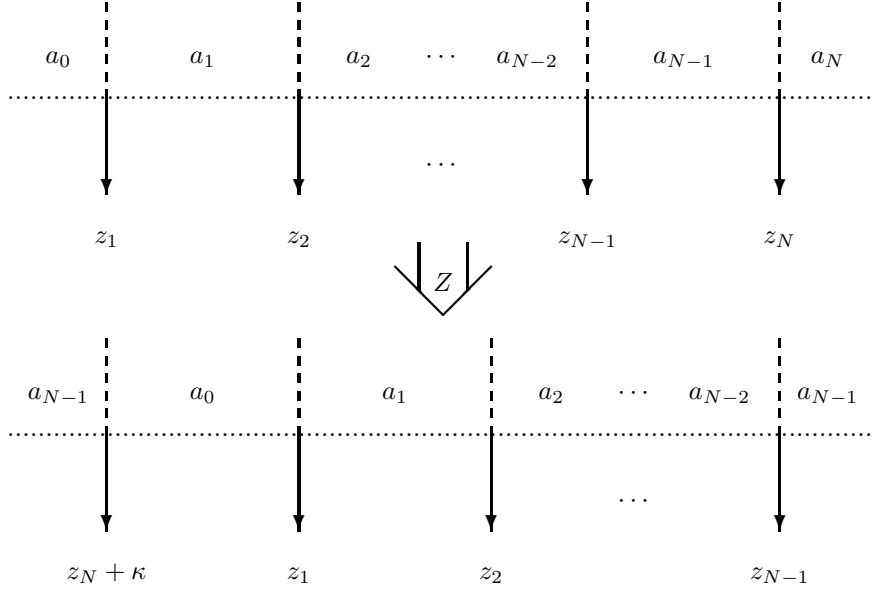

\begin{defn}
A {\em boundary operator} 
$Z = Z^{l_1,\ldots,l_N}_{z_1,\ldots,z_N}(\kappa,c)$ is a linear map from 
$\Bethe^{l_1,\ldots,l_N}_{z_1,\ldots,z_N}$ to
$\Bethe^{l_N, l_1,\ldots,l_{N-1}}_{z_N, z_1,\ldots,z_{N-1}}$ defined by
\begin{multline}
    Z (
    e^{2\pi i a_N \nu \eta } |a_0, \ldots, a_N; z_1, \ldots, z_N \rangle
    ) = \\
    =
    e^{2\pi i a_{N-1} \nu \eta } e^{c(a_N - a_{N-1} - 2 l_N)}
    |a_{N-1}, a_0, \ldots, a_{N-1};
     z_N + \kappa, z_1, \ldots, z_{N-1} \rangle.
\label{boundary-on-bethe}
\end{multline}
(Note that $a_0 = a_N$ by \defref{def:bethe}. cf.~\figref{fig:Z})
\label{def:boundary}
\end{defn}

With respect to the basis \eqref{def:bethe-basis}, $Z$ is a composite of a
permutation and a diagonal matrix.

The following property is important.
\begin{lem}
The operator $Z$ commutes with an $R$ matrix as follows:
\begin{multline}
    \check R^{l_{N-1},l_N}_{12}(z_{N-1} - z_N)
    Z^{l_N,\       l_1, \ldots, l_{N-1}}_{%
       z_N+\kappa, z_1, \ldots, z_{N-1}}
      (\kappa, c)
    Z^{l_1, \ldots, l_N}_{z_1, \ldots, z_N} (\kappa, c) = \\
    =
    Z^{l_{N-1},\       l_1, \ldots, l_{N-2}, l_N}_{%
       z_{N-1}+\kappa, z_1, \ldots, z_{N-2}, z_N} (\kappa, c)
    Z^{l_1, \ldots, l_{N-2}, l_N, l_{N-1}}_{%
       z_1, \ldots, z_{N-2}, z_N, z_{N-1}} (\kappa, c)
    \check R^{l_{N-1},l_N}_{N-1,N}(z_{N-1} - z_N)
\end{multline}
as a map from $\Bethe^{l_1, \ldots, l_N}_{z_1, \ldots, z_N}$ to 
$\Bethe^{l_N, l_{N-1}, l_1, \ldots, l_{N-2}}_{%
         z_N, z_{N-1}, z_1, \ldots, z_{N-2}}$.
\label{commutativity-of-Z-R}
\end{lem}
This is easily proved by comparing the action of the both hand sides on
the basis \eqref{def:bethe-basis}. A graphical notation for this
equality is shown in \figref{fig:RZZ=ZZR}.

\begin{figure}[t]
\setlength{\unitlength}{0.00083300in}%
\begingroup\makeatletter\ifx\SetFigFont\undefined
% extract first six characters in \fmtname
\def\x#1#2#3#4#5#6#7\relax{\def\x{#1#2#3#4#5#6}}%
\expandafter\x\fmtname xxxxxx\relax \def\y{splain}%
\ifx\x\y   % LaTeX or SliTeX?
\gdef\SetFigFont#1#2#3{%
  \ifnum #1<17\tiny\else \ifnum #1<20\small\else
  \ifnum #1<24\normalsize\else \ifnum #1<29\large\else
  \ifnum #1<34\Large\else \ifnum #1<41\LARGE\else
     \huge\fi\fi\fi\fi\fi\fi
  \csname #3\endcsname}%
\else
\gdef\SetFigFont#1#2#3{\begingroup
  \count@#1\relax \ifnum 25<\count@\count@25\fi
  \def\x{\endgroup\@setsize\SetFigFont{#2pt}}%
  \expandafter\x
    \csname \romannumeral\the\count@ pt\expandafter\endcsname
    \csname @\romannumeral\the\count@ pt\endcsname
  \csname #3\endcsname}%
\fi
\fi\endgroup
\begin{picture}(5716,4224)(2243,-3973)
\thicklines
\put(7801,-661){\circle{300}}
\put(2401,-661){\circle{300}}
\put(2401,-1261){\circle{300}}
\put(7801,-1261){\circle{300}}
\put(2401,-3061){\circle{300}}
\put(2401,-3661){\circle{300}}
\put(7801,-3061){\circle{300}}
\put(7801,-3661){\circle{300}}
\put(7201,-361){\line( 0,-1){300}}
\put(7201,-661){\vector( 1, 0){450}}
\multiput(7201,-361)(0.00000,109.09091){6}{\line( 0, 1){ 54.545}}
\put(6601,-361){\line( 0,-1){900}}
\put(6601,-1261){\vector( 1, 0){1050}}
\multiput(6601,-361)(0.00000,109.09091){6}{\line( 0, 1){ 54.545}}
\put(4201,-361){\vector( 0,-1){1200}}
\multiput(4201,-361)(0.00000,109.09091){6}{\line( 0, 1){ 54.545}}
\put(7201,-2761){\line( 0,-1){900}}
\put(7201,-3661){\vector( 1, 0){450}}
\multiput(7201,-2761)(0.00000,109.09091){6}{\line( 0, 1){ 54.545}}
\put(6601,-2761){\line( 0,-1){300}}
\put(6601,-3061){\vector( 1, 0){1050}}
\multiput(6601,-2761)(0.00000,109.09091){6}{\line( 0, 1){ 54.545}}
\put(4201,-2761){\vector( 0,-1){1200}}
\multiput(4201,-2761)(0.00000,109.09091){6}{\line( 0, 1){ 54.545}}
\put(2551,-661){\line( 1, 0){450}}
\put(3001,-661){\vector( 0,-1){900}}
\put(2551,-1261){\line( 1, 0){1050}}
\put(3601,-1261){\vector( 0,-1){300}}
\put(2551,-3061){\line( 1, 0){1050}}
\put(3601,-3061){\vector( 0,-1){900}}
\put(2551,-3661){\line( 1, 0){450}}
\put(3001,-3661){\vector( 0,-1){300}}
\put(5026,-1636){\line( 0,-1){225}}
\put(5101,-1636){\line( 0,-1){225}}
\multiput(2401,-361)(47.78761,0.00000){114}{\makebox(6.6667,10.0000){\SetFigFont{10}{12}{rm}.}}
\multiput(2401,-2761)(47.78761,0.00000){114}{\makebox(6.6667,10.0000){\SetFigFont{10}{12}{rm}.}}
\put(7501,-136){\makebox(0,0)[b]{\smash{$a_N$}}}
\put(6901,-136){\makebox(0,0)[b]{\smash{$a_{N-1}$}}}
\put(6301,-136){\makebox(0,0)[b]{\smash{$a_{N-2}$}}}
\put(5401,-136){\makebox(0,0)[b]{\smash{$\cdots$}}}
\put(5401,-811){\makebox(0,0)[b]{\smash{$\cdots$}}}
\put(7501,-2536){\makebox(0,0)[b]{\smash{$a_N$}}}
\put(6901,-2536){\makebox(0,0)[b]{\smash{$a_{N-1}$}}}
\put(6301,-2536){\makebox(0,0)[b]{\smash{$a_{N-2}$}}}
\put(5401,-2536){\makebox(0,0)[b]{\smash{$\cdots$}}}
\put(5401,-3211){\makebox(0,0)[b]{\smash{$\cdots$}}}
\put(7801,-736){\makebox(0,0)[b]{\smash{$Z$}}}
\put(2401,-736){\makebox(0,0)[b]{\smash{$Z$}}}
\put(4351,-3361){\makebox(0,0)[b]{\smash{$z_1$}}}
\put(7426,-811){\makebox(0,0)[b]{\smash{$z_N$}}}
\put(2401,-1336){\makebox(0,0)[b]{\smash{$Z$}}}
\put(7801,-1336){\makebox(0,0)[b]{\smash{$Z$}}}
\put(4351,-961){\makebox(0,0)[b]{\smash{$z_1$}}}
\put(2851,-1411){\makebox(0,0)[b]{\smash{$R$}}}
\put(2401,-3136){\makebox(0,0)[b]{\smash{$Z$}}}
\put(2401,-3736){\makebox(0,0)[b]{\smash{$Z$}}}
\put(7801,-3136){\makebox(0,0)[b]{\smash{$Z$}}}
\put(7801,-3736){\makebox(0,0)[b]{\smash{$Z$}}}
\put(7051,-3436){\makebox(0,0)[b]{\smash{$z_N$}}}
\put(7426,-3211){\makebox(0,0)[b]{\smash{$z_{N-1}$}}}
\put(7051,-2986){\makebox(0,0)[b]{\smash{$R$}}}
\put(3076,-3211){\makebox(0,0)[b]{\smash{$z_{N-1}+\kappa$}}}
\put(7201,-1186){\makebox(0,0)[b]{\smash{$z_{N-1}$}}}
\put(2851,-3586){\makebox(0,0)[b]{\smash{$z_N+\kappa$}}}
\put(2851,-586){\makebox(0,0)[b]{\smash{$z_N+\kappa$}}}
\put(3451,-1186){\makebox(0,0)[b]{\smash{$z_{N-1}+\kappa$}}}
\put(3376,-2536){\makebox(0,0)[b]{\smash{$a_0$}}}
\put(4501,-2536){\makebox(0,0)[b]{\smash{$a_1$}}}
\put(4501,-136){\makebox(0,0)[b]{\smash{$a_1$}}}
\put(3376,-136){\makebox(0,0)[b]{\smash{$a_0$}}}
\end{picture}
\caption{$R Z Z = Z Z R$, \lemref{commutativity-of-Z-R}.}
\label{fig:RZZ=ZZR}
\end{figure}

\section{Difference equation}
\setcounter{equation}{0}

In this section we introduce a system of difference equations which is
an elliptic analogue of the equations introduced in \cite{fre-resh},
and show holonomicity of this system in the sense of
Aomoto \cite{aom:88}.

Let us first define a linear operator $A_j(\vec z)$, 
$\vec z = (z_1, \ldots, z_N)$, from 
$\Bethe^{l_1, \ldots, l_N}_{z_1, \ldots, z_N}$ to
$\Bethe^{l_1, \ldots, l_j,        \ldots, l_N}_{%
         z_1, \ldots, z_j+\kappa, \ldots, z_N}$ by
\begin{multline}
    A_j(\vec z) := 
    \check R^{l_j, l_{j-1}}_{j-1,j}(z_j + \kappa - z_{j-1})
    \ldots
    \check R^{l_j, l_1}_{12}(z_j + \kappa - z_1) \cdot \\
    \cdot
    Z^{l_1, \ldots, l_{j-1}, l_{j+1}, \ldots, l_N, l_j}_{%
       z_1, \ldots, z_{j-1}, z_{j+1}, \ldots, z_N, z_j} (\kappa, c)
    \check R^{l_j, l_N}_{N-1,N}(z_j - z_N)
    \ldots
    \check R^{l_j, l_{j+1}}_{j,j+1}(z_j - z_{j+1}),
\label{def:a-j:R-check}
\end{multline}
or, equivalently,
\begin{multline}
    A_j(\vec z) := 
    R^{l_j, l_{j-1}}_{j,j-1}(z_j + \kappa - z_{j-1})
    \ldots
    R^{l_j, l_1}_{j,1}(z_j + \kappa - z_1) \\
    Z^{(j)} (\kappa, c) 
    R^{l_j, l_N}_{j,N}(z_j - z_N)
    \ldots
    R^{l_j, l_{j+1}}_{j,j+1}(z_j - z_{j+1}),
\label{def:a-j:R}
\end{multline}
where
\begin{multline}
    Z^{(j)} (\kappa, c) =
    P_{j-1, j} P_{j-2, j-1} \cdots P_{12} \\
    Z^{l_1, \ldots, l_{j-1}, l_{j+1}, \ldots, l_N, l_j}_{%
       z_1, \ldots, z_{j-1}, z_{j+1}, \ldots, z_N, z_j} (\kappa, c)
    P_{N-1, N} \cdots P_{j+1, j+2} P_{j, j+1}.
\label{def:Z-j}
\end{multline}

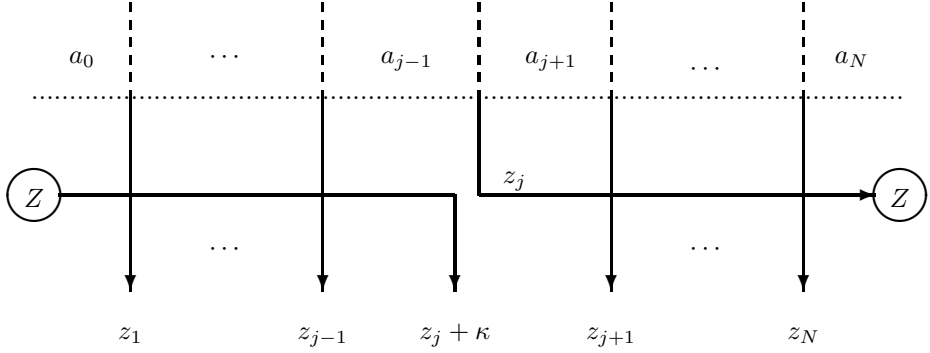
\begin{figure}[ht]
\setlength{\unitlength}{0.00083300in}%
\begingroup\makeatletter\ifx\SetFigFont\undefined
% extract first six characters in \fmtname
\def\x#1#2#3#4#5#6#7\relax{\def\x{#1#2#3#4#5#6}}%
\expandafter\x\fmtname xxxxxx\relax \def\y{splain}%
\ifx\x\y   % LaTeX or SliTeX?
\gdef\SetFigFont#1#2#3{%
  \ifnum #1<17\tiny\else \ifnum #1<20\small\else
  \ifnum #1<24\normalsize\else \ifnum #1<29\large\else
  \ifnum #1<34\Large\else \ifnum #1<41\LARGE\else
     \huge\fi\fi\fi\fi\fi\fi
  \csname #3\endcsname}%
\else
\gdef\SetFigFont#1#2#3{\begingroup
  \count@#1\relax \ifnum 25<\count@\count@25\fi
  \def\x{\endgroup\@setsize\SetFigFont{#2pt}}%
  \expandafter\x
    \csname \romannumeral\the\count@ pt\expandafter\endcsname
    \csname @\romannumeral\the\count@ pt\endcsname
  \csname #3\endcsname}%
\fi
\fi\endgroup
\begin{picture}(5716,2148)(2243,-1897)
\thicklines
\put(2401,-961){\circle{300}}
\put(7801,-961){\circle{300}}
\multiput(3001,-361)(0.00000,109.09091){6}{\line( 0, 1){ 54.545}}
\put(3001,-361){\vector( 0,-1){1200}}
\put(2551,-961){\line( 1, 0){2475}}
\put(5026,-961){\vector( 0,-1){600}}
\multiput(5176,-361)(0.00000,109.09091){6}{\line( 0, 1){ 54.545}}
\put(5176,-361){\line( 0,-1){600}}
\put(5176,-961){\vector( 1, 0){2475}}
\multiput(4201,-361)(0.00000,109.09091){6}{\line( 0, 1){ 54.545}}
\put(4201,-361){\vector( 0,-1){1200}}
\multiput(6001,-361)(0.00000,109.09091){6}{\line( 0, 1){ 54.545}}
\put(6001,-361){\vector( 0,-1){1200}}
\multiput(7201,-361)(0.00000,109.09091){6}{\line( 0, 1){ 54.545}}
\put(7201,-361){\vector( 0,-1){1200}}
\multiput(2401,-361)(47.78761,0.00000){114}{\makebox(6.6667,10.0000){\SetFigFont{10}{12}{rm}.}}
\put(2401,-1036){\makebox(0,0)[b]{\smash{$Z$}}}
\put(3001,-1861){\makebox(0,0)[b]{\smash{$z_1$}}}
\put(6001,-1861){\makebox(0,0)[b]{\smash{$z_{j+1}$}}}
\put(7201,-1861){\makebox(0,0)[b]{\smash{$z_N$}}}
\put(7801,-1036){\makebox(0,0)[b]{\smash{$Z$}}}
\put(3601,-136){\makebox(0,0)[b]{\smash{$\cdots$}}}
\put(6601,-211){\makebox(0,0)[b]{\smash{$\cdots$}}}
\put(3601,-1336){\makebox(0,0)[b]{\smash{$\cdots$}}}
\put(6601,-1336){\makebox(0,0)[b]{\smash{$\cdots$}}}
\put(4201,-1861){\makebox(0,0)[b]{\smash{$z_{j-1}$}}}
\put(5026,-1861){\makebox(0,0)[b]{\smash{$z_j+\kappa$}}}
\put(5326,-886){\makebox(0,0)[lb]{\smash{$z_j$}}}
\put(2701,-136){\makebox(0,0)[b]{\smash{$a_0$}}}
\put(4726,-136){\makebox(0,0)[b]{\smash{$a_{j-1}$}}}
\put(5626,-136){\makebox(0,0)[b]{\smash{$a_{j+1}$}}}
\put(7501,-136){\makebox(0,0)[b]{\smash{$a_N$}}}
\end{picture}
\caption{Operator $A_j(\vec z)$}
\label{fig:A-j}
\end{figure}

Then we can define our main object, a system of difference equations:
\begin{equation}
    f({\vec z}_j)
    =
    A_j(\vec z) f(\vec z),
\label{diff-eq}
\end{equation}
for $j = 1, \ldots, N$ where 
${\vec z}_j = (z_1, \ldots, z_j + \kappa, \ldots, z_N)$ and
$f(\vec z) \in \Bethe^{l_1, \ldots, l_N}_{z_1, \ldots, z_N}$.
Expanding $f$ as in \eqref{def:bethe},
\begin{equation*}
    f(\vec z) = \sum_{a=0}^{r-1} e^{2\pi i a \nu \eta}
             \sum_{\vec a} f_{\vec a} (\vec z) \,
             | \vec a; \vec z \rangle, 
\end{equation*}
we can regard the system \eqref{diff-eq} as a system of equations for
$f_{\vec a}(\vec z)$. In this IRF picture, this system is holonomic in the
sense of Aomoto \cite{aom:88} due to the following:

\begin{prop}
\label{holonomic}
Operators $A_j(\vec z)$ are compatible:
$$
    A_j({\vec z}_k) A_k(\vec z) = A_k({\vec z}_j) A_j (\vec z).
$$
\end{prop}

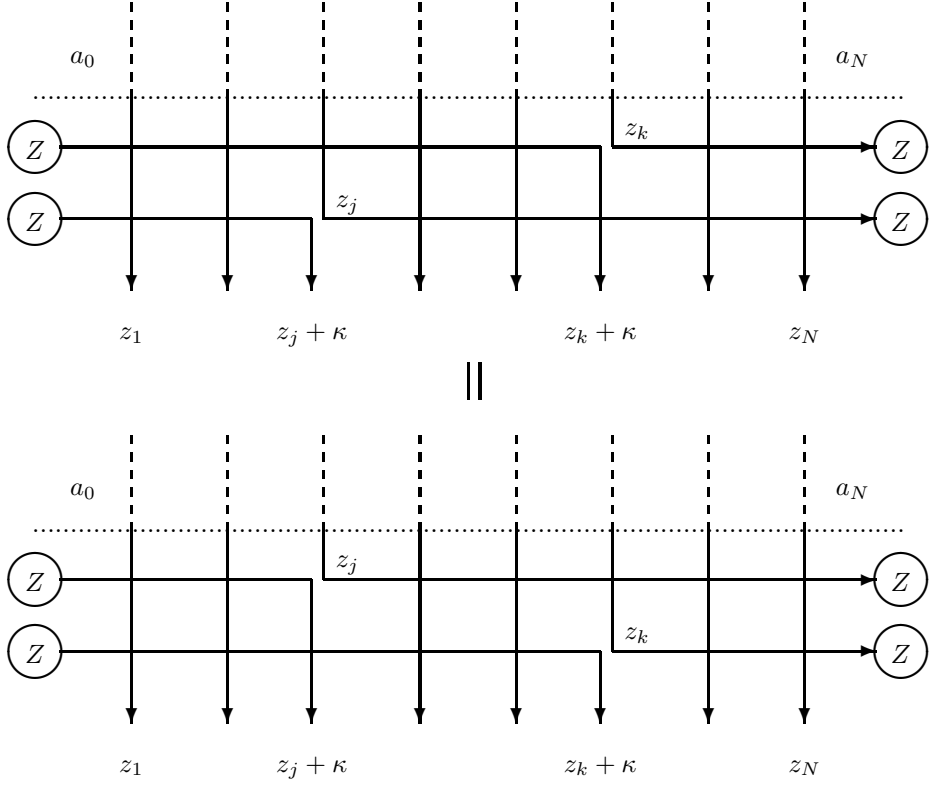
\begin{figure}[ht]
\setlength{\unitlength}{0.00083300in}%
\begingroup\makeatletter\ifx\SetFigFont\undefined
% extract first six characters in \fmtname
\def\x#1#2#3#4#5#6#7\relax{\def\x{#1#2#3#4#5#6}}%
\expandafter\x\fmtname xxxxxx\relax \def\y{splain}%
\ifx\x\y   % LaTeX or SliTeX?
\gdef\SetFigFont#1#2#3{%
  \ifnum #1<17\tiny\else \ifnum #1<20\small\else
  \ifnum #1<24\normalsize\else \ifnum #1<29\large\else
  \ifnum #1<34\Large\else \ifnum #1<41\LARGE\else
     \huge\fi\fi\fi\fi\fi\fi
  \csname #3\endcsname}%
\else
\gdef\SetFigFont#1#2#3{\begingroup
  \count@#1\relax \ifnum 25<\count@\count@25\fi
  \def\x{\endgroup\@setsize\SetFigFont{#2pt}}%
  \expandafter\x
    \csname \romannumeral\the\count@ pt\expandafter\endcsname
    \csname @\romannumeral\the\count@ pt\endcsname
  \csname #3\endcsname}%
\fi
\fi\endgroup
\begin{picture}(5716,4848)(2243,-4597)
\thicklines
\put(2401,-661){\circle{300}}
\put(7801,-661){\circle{300}}
\put(7801,-1111){\circle{300}}
\put(2401,-1111){\circle{300}}
\put(7801,-3361){\circle{300}}
\put(7801,-3811){\circle{300}}
\put(2401,-3361){\circle{300}}
\put(2401,-3811){\circle{300}}
\multiput(3001,-361)(0.00000,109.09091){6}{\line( 0, 1){ 54.545}}
\put(3001,-361){\vector( 0,-1){1200}}
\multiput(7201,-361)(0.00000,109.09091){6}{\line( 0, 1){ 54.545}}
\put(7201,-361){\vector( 0,-1){1200}}
\multiput(6601,-361)(0.00000,109.09091){6}{\line( 0, 1){ 54.545}}
\put(2551,-1111){\line( 1, 0){1575}}
\put(4126,-1111){\vector( 0,-1){450}}
\put(4201,-361){\line( 0,-1){750}}
\put(4201,-1111){\vector( 1, 0){3450}}
\multiput(6001,-361)(0.00000,109.09091){6}{\line( 0, 1){ 54.545}}
\put(2551,-661){\line( 1, 0){3375}}
\put(5926,-661){\vector( 0,-1){900}}
\put(6001,-361){\line( 0,-1){300}}
\put(6001,-661){\vector( 1, 0){1650}}
\put(6601,-361){\vector( 0,-1){1200}}
\put(5401,-361){\vector( 0,-1){1200}}
\multiput(5401,-361)(0.00000,109.09091){6}{\line( 0, 1){ 54.545}}
\multiput(4201,-361)(0.00000,109.09091){6}{\line( 0, 1){ 54.545}}
\put(3601,-361){\vector( 0,-1){1200}}
\multiput(3601,-361)(0.00000,109.09091){6}{\line( 0, 1){ 54.545}}
\put(4801,-361){\vector( 0,-1){1200}}
\multiput(4801,-361)(0.00000,109.09091){6}{\line( 0, 1){ 54.545}}
\multiput(3001,-3061)(0.00000,109.09091){6}{\line( 0, 1){ 54.545}}
\put(3001,-3061){\vector( 0,-1){1200}}
\multiput(7201,-3061)(0.00000,109.09091){6}{\line( 0, 1){ 54.545}}
\put(7201,-3061){\vector( 0,-1){1200}}
\multiput(6601,-3061)(0.00000,109.09091){6}{\line( 0, 1){ 54.545}}
\put(2551,-3361){\line( 1, 0){1575}}
\put(4126,-3361){\vector( 0,-1){900}}
\put(4201,-3061){\line( 0,-1){300}}
\put(4201,-3361){\vector( 1, 0){3450}}
\multiput(6001,-3061)(0.00000,109.09091){6}{\line( 0, 1){ 54.545}}
\put(6001,-3061){\line( 0,-1){750}}
\put(6001,-3811){\vector( 1, 0){1650}}
\put(6601,-3061){\vector( 0,-1){1200}}
\put(5401,-3061){\vector( 0,-1){1200}}
\multiput(5401,-3061)(0.00000,109.09091){6}{\line( 0, 1){ 54.545}}
\multiput(4201,-3061)(0.00000,109.09091){6}{\line( 0, 1){ 54.545}}
\put(3601,-3061){\vector( 0,-1){1200}}
\multiput(3601,-3061)(0.00000,109.09091){6}{\line( 0, 1){ 54.545}}
\put(4801,-3061){\vector( 0,-1){1200}}
\multiput(4801,-3061)(0.00000,109.09091){6}{\line( 0, 1){ 54.545}}
\put(5101,-2011){\line( 0,-1){225}}
\put(5176,-2011){\line( 0,-1){225}}
\put(2551,-3811){\line( 1, 0){3375}}
\put(5926,-3811){\vector( 0,-1){450}}
\multiput(2401,-361)(47.78761,0.00000){114}{\makebox(6.6667,10.0000){\SetFigFont{10}{12}{rm}.}}
\multiput(2401,-3061)(47.78761,0.00000){114}{\makebox(6.6667,10.0000){\SetFigFont{10}{12}{rm}.}}
\put(2701,-136){\makebox(0,0)[b]{\smash{$a_0$}}}
\put(7501,-136){\makebox(0,0)[b]{\smash{$a_N$}}}
\put(2401,-736){\makebox(0,0)[b]{\smash{$Z$}}}
\put(7801,-736){\makebox(0,0)[b]{\smash{$Z$}}}
\put(7801,-1186){\makebox(0,0)[b]{\smash{$Z$}}}
\put(2401,-1186){\makebox(0,0)[b]{\smash{$Z$}}}
\put(3001,-1861){\makebox(0,0)[b]{\smash{$z_1$}}}
\put(7201,-1861){\makebox(0,0)[b]{\smash{$z_N$}}}
\put(4126,-1861){\makebox(0,0)[b]{\smash{$z_j+\kappa$}}}
\put(4276,-1036){\makebox(0,0)[lb]{\smash{$z_j$}}}
\put(6076,-586){\makebox(0,0)[lb]{\smash{$z_k$}}}
\put(5926,-1861){\makebox(0,0)[b]{\smash{$z_k+\kappa$}}}
\put(2701,-2836){\makebox(0,0)[b]{\smash{$a_0$}}}
\put(7501,-2836){\makebox(0,0)[b]{\smash{$a_N$}}}
\put(3001,-4561){\makebox(0,0)[b]{\smash{$z_1$}}}
\put(7201,-4561){\makebox(0,0)[b]{\smash{$z_N$}}}
\put(4126,-4561){\makebox(0,0)[b]{\smash{$z_j+\kappa$}}}
\put(7801,-3436){\makebox(0,0)[b]{\smash{$Z$}}}
\put(7801,-3886){\makebox(0,0)[b]{\smash{$Z$}}}
\put(2401,-3436){\makebox(0,0)[b]{\smash{$Z$}}}
\put(2401,-3886){\makebox(0,0)[b]{\smash{$Z$}}}
\put(5926,-4561){\makebox(0,0)[b]{\smash{$z_k+\kappa$}}}
\put(6076,-3736){\makebox(0,0)[lb]{\smash{$z_k$}}}
\put(4276,-3286){\makebox(0,0)[lb]{\smash{$z_j$}}}
\end{picture}
\caption{Holonomicity of the system \eqref{diff-eq}.}
\label{fig:holon}
\end{figure}

This follows from the Yang-Baxter equation \eqref{yb}, the unitarity
\eqref{unitarity} and the commutativity of $R$ and $Z$
(\lemref{commutativity-of-Z-R}) as in Theorem 5.4 of
\cite{fre-resh}. Symbolically, we have only to change the order of
crossings of lines in \figref{fig:holon}, using the procedures
\figref{fig:yb}, \figref{fig:unitarity} and \figref{fig:RZZ=ZZR}.

\section{Jackson-type integral solution}
\setcounter{equation}{0}

In this section we give a solution of the system of difference equations
\eqref{diff-eq} which is expressed as a Jackson-type integral of Bethe
vectors. This result is an elliptic analogue of that of \cite{resh:92}.
We assume that $M := l_1 + \cdots + l_N$ is an integer.

\begin{figure}[ht]
\setlength{\unitlength}{0.00083300in}%
\begingroup\makeatletter\ifx\SetFigFont\undefined
% extract first six characters in \fmtname
\def\x#1#2#3#4#5#6#7\relax{\def\x{#1#2#3#4#5#6}}%
\expandafter\x\fmtname xxxxxx\relax \def\y{splain}%
\ifx\x\y   % LaTeX or SliTeX?
\gdef\SetFigFont#1#2#3{%
  \ifnum #1<17\tiny\else \ifnum #1<20\small\else
  \ifnum #1<24\normalsize\else \ifnum #1<29\large\else
  \ifnum #1<34\Large\else \ifnum #1<41\LARGE\else
     \huge\fi\fi\fi\fi\fi\fi
  \csname #3\endcsname}%
\else
\gdef\SetFigFont#1#2#3{\begingroup
  \count@#1\relax \ifnum 25<\count@\count@25\fi
  \def\x{\endgroup\@setsize\SetFigFont{#2pt}}%
  \expandafter\x
    \csname \romannumeral\the\count@ pt\expandafter\endcsname
    \csname @\romannumeral\the\count@ pt\endcsname
  \csname #3\endcsname}%
\fi
\fi\endgroup
\begin{picture}(3757,1473)(3656,-1522)
\thicklines
\multiput(5101,-661)(-120.00000,0.00000){3}{\line(-1, 0){ 60.000}}
\put(5701,-61){\vector( 0,-1){1200}}
\multiput(5101,-1261)(0.00000,48.00000){26}{\makebox(6.6667,10.0000){\SetFigFont{10}{12}{rm}.}}
\put(5101,-661){\vector( 1, 0){3000}}
\put(6301,-61){\vector( 0,-1){1200}}
\multiput(8401,-661)(-120.00000,0.00000){3}{\line(-1, 0){ 60.000}}
\put(7501,-61){\vector( 0,-1){1200}}
\multiput(8101,-1261)(0.00000,48.00000){26}{\makebox(6.6667,10.0000){\SetFigFont{10}{12}{rm}.}}
\put(5026,-1036){\makebox(0,0)[rb]{\smash{$a'+1$}}}
\put(5026,-361){\makebox(0,0)[rb]{\smash{$a'$}}}
\put(8176,-361){\makebox(0,0)[lb]{\smash{$a$}}}
\put(8176,-1036){\makebox(0,0)[lb]{\smash{$a-1$}}}
\put(6301,-1486){\makebox(0,0)[b]{\smash{$z_j$}}}
\put(5701,-1486){\makebox(0,0)[b]{\smash{$z_1$}}}
\put(7501,-1486){\makebox(0,0)[b]{\smash{$z_N$}}}
\put(6901,-1036){\makebox(0,0)[b]{\smash{$\cdots$}}}
\put(6901,-436){\makebox(0,0)[b]{\smash{$\cdots$}}}
\put(5401,-586){\makebox(0,0)[b]{\smash{$u$}}}
\end{picture}
\caption{Operator $B^{l_N, \ldots, l_1}_{a,a'}(u; z_N, \ldots, z_1)$.}
\label{fig:B-a-a'}
\end{figure}
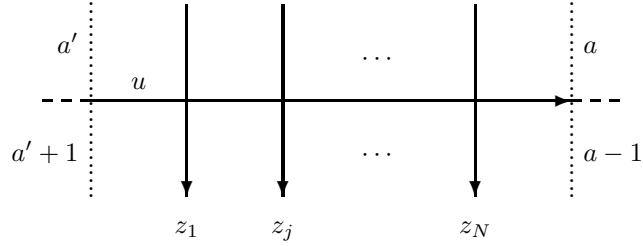

First let us define the {\em monodromy matrix} $T(u;z_1, \ldots, z_N)$ by
\begin{multline}
    T^{l_N, \ldots, l_1}(u;z_N, \ldots, z_1) = \\
    =
    \begin{pmatrix}
    A^{l_N, \ldots, l_1}(u;z_N, \ldots, z_1) &
    B^{l_N, \ldots, l_1}(u;z_N, \ldots, z_1) \\
    C^{l_N, \ldots, l_1}(u;z_N, \ldots, z_1) & 
    D^{l_N, \ldots, l_1}(u;z_N, \ldots, z_1) 
    \end{pmatrix}
    := \\
    :=
    L_N^{l_N}(u-z_N) \cdots L_2^{l_2}(u-z_2) L_1^{l_1}(u-z_1)
\label{def:monodromy}
\end{multline}
as an endomorphism of
$ \Complex^2 \tensor \calH 
= \Complex^2 \tensor V^{l_1} \tensor \cdots \tensor V^{l_N}$.
Here $L_j^{l_j}(u)$ acts non-trivially only on $\Complex^2$ and
$V^{l_j}$:
\begin{gather*}
    L_j^{l_j}(u) = 
    \sum_{a=0}^3 W_a^L(u) \sigma^a \tensor \rho_j^{l_j}(S^a),
\\
    \rho_j^{l_j} = 1 \tensor \cdots \tensor 1 \tensor \rho^{l_j}
           \tensor 1 \tensor \cdots \tensor 1:
    U_{\tau,\eta}(sl(2)) \to \End(\calH).
\end{gather*}
Bethe vectors of the XYZ type spin chains are constructed by means of the
twisted monodromy matrix defined as follows:
\begin{multline}
    T^{l_N, \ldots, l_1}_{\lam,\lam'}(u;z_N, \ldots, z_1) = \\
    =
    \begin{pmatrix}
    A^{l_N, \ldots, l_1}_{\lam,\lam'}(u;z_N, \ldots, z_1) &
    B^{l_N, \ldots, l_1}_{\lam,\lam'}(u;z_N, \ldots, z_1) \\
    C^{l_N, \ldots, l_1}_{\lam,\lam'}(u;z_N, \ldots, z_1) &  
    D^{l_N, \ldots, l_1}_{\lam,\lam'}(u;z_N, \ldots, z_1)
    \end{pmatrix} := \\
    :=
    M_\lam(u)^{-1} T^{l_N, \ldots, l_1}(u;z_N, \ldots, z_1) M_{\lam'}(u),
\label{def:twisted-monodromy}
\end{multline}
where $M_\lam(u)$ is defined by \eqref{def:gauge-trans}. We often
denote, for example, 
$B^{l_N, \ldots, l_1}_{\lamcirc+ 2a\eta, \lamcirc + 2a'\eta}$
by $B^{l_N, \ldots, l_1}_{a,a'}$ or $B_{a,a'}$ for simplicity, and
graphically as in \figref{fig:B-a-a'}. 
(cf.\ \figref{fig:twisted-L-on-int-vec}.)
It follows from \eqref{alpha-on-vac}, \eqref{gamma-on-vac},
\eqref{delta-on-vac} that
\begin{align}
    A^{l_N, \ldots, l_1}_{a + 2M,a}(u;z_N, \ldots, z_1) \Omega_a(\vec z)
    &= \left( \prod_{j=1}^N \alpha^{l_j}(u - z_j) \right)
       \Omega_{a-1}(\vec z),
\label{A-on-vac}
\\
    C^{l_N, \ldots, l_1}_{a + 2M,a}(u;z_N, \ldots, z_1) \Omega_a(\vec z)
    &= 0
\label{C-on-vac}
\\
    D^{l_N, \ldots, l_1}_{a + 2M,a}(u;z_N, \ldots, z_1) \Omega_a(\vec z)
    &= \left( \prod_{j=1}^N \delta^{l_j}(u - z_j) \right)
       \Omega_{a+1}(\vec z).
\label{D-on-vac}
\end{align}
By a standard argument in \cite{takh-fad}, these operators satisfy the
following commutation relations:
\begin{gather}
\begin{split}
    B_{a,a'+1}&(u;z_N, \ldots, z_1) B_{a+1, a'}(v;z_N, \ldots, z_1)
    = \\
    = &
    B_{a,a'+1}(v;z_N, \ldots, z_1) B_{a+1, a'}(u;z_N, \ldots, z_1)
\end{split}
\label{comm-rel-B-B}
\\
\begin{split}
    A_{a,a'+1}&(u;z_N, \ldots, z_1) B_{a+1, a'}(v;z_N, \ldots, z_1)
    = \\
    =& \alpha(u-v)
    B_{a,a'-1}(v;z_N, \ldots, z_1) A_{a+1,a'}(u;z_N, \ldots, z_1)
    +\\
    +& \beta_{a'}(u-v)
    B_{a,a'-1}(u;z_N, \ldots, z_1) A_{a+1,a'}(v;z_N, \ldots, z_1),
\end{split}
\label{comm-rel-A-B}
\\
\begin{split}
    D_{a-1,a'}&(u;z_N, \ldots, z_1) B_{a, a'-1}(v;z_N, \ldots, z_1)
    = \\
    =& \alpha(v-u)
    B_{a+1,a'}(v;z_N, \ldots, z_1) D_{a,a'-1}(u;z_N, \ldots, z_1)
    -\\
    -& \beta_{a }(u-v)
    B_{a+1,a'}(u;z_N, \ldots, z_1) D_{a,a'-1}(v;z_N, \ldots, z_1),
\end{split}
\label{comm-rel-D-B}
\end{gather}
where
\begin{equation}
    \alpha(u) = \frac{\theta_{11}(u - 2\eta)}{\theta_{11}(u)},
    \qquad
    \beta_a(u) = 
    \frac{\theta_{11}(u - \lamcirc - 2 a \eta) \theta_{11}(2\eta)}
         {\theta_{11}(u) \theta_{11}(\lamcirc + 2 a \eta)}.
\end{equation}

Let us recall the definition of $N$-cycle with step $\kappa$ in
\cite{resh:92}.

\begin{defn}
Let $\Psi$ be a function on $\Complex^N$ and $\calC$ be a subset of 
$\{ (x_1 + m_1 \kappa, \ldots, x_N + m_N \kappa ) \, |\,
    (m_1, \ldots, m_N)\in \Integer^N \}$ for a certain
$(x_1, \ldots, x_N)$. We call $\calC$
an {\em $N$-cycle with step $\kappa$ for $\Psi$} if
\begin{equation*}
    \sum_{\vec t \in \calC} \Psi(\vec t + \vec n \kappa)
    =
    \sum_{\vec t \in \calC} \Psi(\vec t \,),
\end{equation*}
for any $\vec n \in \Integer^N$.
\end{defn}

Functions $F^l(t)$ and $\Phi(t)$ are defined as solutions of the
ordinary difference equations:
\begin{align}
    F^l(t+\kappa)  &= \frac{\delta^l(t)}{\alpha^l(t+\kappa)} F^l(t), 
\label{diff-eq:F-l}
\\
    \Phi(t+\kappa) &= \frac{\alpha(-t)}{\alpha(t+\kappa)} \Phi(t).
\label{diff-eq:Phi}
\end{align}
We shall give an explicit solution to these equations later in terms
of infinite products.

The main theorem of this paper is:
\begin{thm}
Let
\begin{equation}
    \varphi(\vec z | \vec t \,) = 
    e^{c \sum_{i=1}^M t_i} \cdot 
    \prod_{1 \leq i<j \leq M} \Phi(t_i - t_j) \cdot
    \prod_{j=1}^M \prod_{n=1}^N F^{l_n}(t_j - z_n),
\label{def:varphi}
\end{equation}
where $\vec t = (t_1, \ldots, t_M)$. Then
\begin{equation}
    f(\vec z) =
    \sum_{\vec t\in\calC} \varphi(\vec z | \vec t \,)
    \Psi(\vec t \,)
\label{solution}
\end{equation}
is a solution of \eqref{diff-eq}, where
\begin{equation}
\begin{split}
    \Psi(\vec t\,) =
    \sum_{a=0}^{r-1} e^{2\pi i a \nu \eta} &
    B_{a+1, a-1}(t_1;z_N, \ldots, z_1) 
    B_{a+2, a-2}(t_2;z_N, \ldots, z_1)
    \cdots \\
    & \cdots
    B_{a+M, a-M}(t_M;z_N, \ldots, z_1) \Omega_{a-M}(\vec z)
\end{split}
\label{def:Psi}
\end{equation}
and $\calC$ is an $N$-cycle with step $\kappa$ for 
$\varphi(\vec z| \vec t \,) \Psi(\vec t \,)$.
\label{main-thm}
\end{thm}

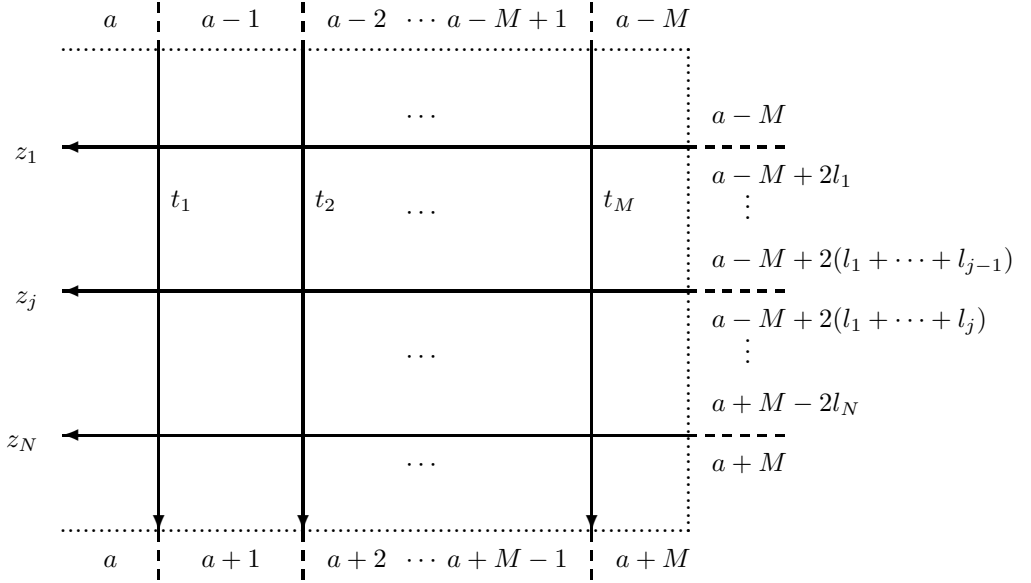
\begin{figure}[ht]
\setlength{\unitlength}{0.00083300in}%
\begingroup\makeatletter\ifx\SetFigFont\undefined
% extract first six characters in \fmtname
\def\x#1#2#3#4#5#6#7\relax{\def\x{#1#2#3#4#5#6}}%
\expandafter\x\fmtname xxxxxx\relax \def\y{splain}%
\ifx\x\y   % LaTeX or SliTeX?
\gdef\SetFigFont#1#2#3{%
  \ifnum #1<17\tiny\else \ifnum #1<20\small\else
  \ifnum #1<24\normalsize\else \ifnum #1<29\large\else
  \ifnum #1<34\Large\else \ifnum #1<41\LARGE\else
     \huge\fi\fi\fi\fi\fi\fi
  \csname #3\endcsname}%
\else
\gdef\SetFigFont#1#2#3{\begingroup
  \count@#1\relax \ifnum 25<\count@\count@25\fi
  \def\x{\endgroup\@setsize\SetFigFont{#2pt}}%
  \expandafter\x
    \csname \romannumeral\the\count@ pt\expandafter\endcsname
    \csname @\romannumeral\the\count@ pt\endcsname
  \csname #3\endcsname}%
\fi
\fi\endgroup
\begin{picture}(4878,3948)(3235,-3997)
\thicklines
\multiput(4201,-361)(0.00000,120.00000){3}{\line( 0, 1){ 60.000}}
\put(4201,-361){\vector( 0,-1){3000}}
\put(7501,-961){\vector(-1, 0){3900}}
\multiput(5101,-361)(0.00000,120.00000){3}{\line( 0, 1){ 60.000}}
\put(5101,-361){\vector( 0,-1){3000}}
\put(6901,-361){\vector( 0,-1){3000}}
\multiput(6901,-361)(0.00000,120.00000){3}{\line( 0, 1){ 60.000}}
\put(7501,-1861){\vector(-1, 0){3900}}
\put(7501,-2761){\vector(-1, 0){3900}}
\multiput(3601,-361)(48.14815,0.00000){82}{\makebox(6.6667,10.0000){\SetFigFont{10}{12}{rm}.}}
\multiput(7501,-361)(0.00000,-47.61905){64}{\makebox(6.6667,10.0000){\SetFigFont{10}{12}{rm}.}}
\multiput(7501,-3361)(-48.14815,0.00000){82}{\makebox(6.6667,10.0000){\SetFigFont{10}{12}{rm}.}}
\multiput(4201,-3661)(0.00000,120.00000){3}{\line( 0, 1){ 60.000}}
\multiput(5101,-3661)(0.00000,120.00000){3}{\line( 0, 1){ 60.000}}
\multiput(6901,-3661)(0.00000,120.00000){3}{\line( 0, 1){ 60.000}}
\multiput(8101,-961)(-109.09091,0.00000){6}{\line(-1, 0){ 54.545}}
\multiput(8101,-1861)(-109.09091,0.00000){6}{\line(-1, 0){ 54.545}}
\multiput(8101,-2761)(-109.09091,0.00000){6}{\line(-1, 0){ 54.545}}
\put(4651,-211){\makebox(0,0)[b]{\smash{$a-1$}}}
\put(7051,-211){\makebox(0,0)[lb]{\smash{$a-M$}}}
\put(5251,-211){\makebox(0,0)[lb]{\smash{$a-2$}}}
\put(6751,-211){\makebox(0,0)[rb]{\smash{$a-M+1$}}}
\put(5851,-211){\makebox(0,0)[b]{\smash{$\cdots$}}}
\put(6976,-1336){\makebox(0,0)[lb]{\smash{$t_M$}}}
\put(5176,-1336){\makebox(0,0)[lb]{\smash{$t_2$}}}
\put(4276,-1336){\makebox(0,0)[lb]{\smash{$t_1$}}}
\put(5851,-811){\makebox(0,0)[b]{\smash{$\cdots$}}}
\put(5851,-1411){\makebox(0,0)[b]{\smash{$\cdots$}}}
\put(5851,-2311){\makebox(0,0)[b]{\smash{$\cdots$}}}
\put(5851,-2986){\makebox(0,0)[b]{\smash{$\cdots$}}}
\put(7651,-811){\makebox(0,0)[lb]{\smash{$a-M$}}}
\put(7651,-1186){\makebox(0,0)[lb]{\smash{$a-M+2l_1$}}}
\put(7651,-1711){\makebox(0,0)[lb]{\smash{$a-M+2(l_1+ \cdots + l_{j-1})$}}}
\put(7876,-1411){\makebox(0,0)[b]{\smash{$\vdots$}}}
\put(7651,-2086){\makebox(0,0)[lb]{\smash{$a-M+2(l_1+ \cdots + l_j)$}}}
\put(7651,-2611){\makebox(0,0)[lb]{\smash{$a+M-2l_N$}}}
\put(7876,-2311){\makebox(0,0)[b]{\smash{$\vdots$}}}
\put(7651,-2986){\makebox(0,0)[lb]{\smash{$a+M$}}}
\put(4651,-3586){\makebox(0,0)[b]{\smash{$a+1$}}}
\put(7051,-3586){\makebox(0,0)[lb]{\smash{$a+M$}}}
\put(6751,-3586){\makebox(0,0)[rb]{\smash{$a+M-1$}}}
\put(5851,-3586){\makebox(0,0)[b]{\smash{$\cdots$}}}
\put(5251,-3586){\makebox(0,0)[lb]{\smash{$a+2$}}}
\put(3901,-211){\makebox(0,0)[b]{\smash{$a$}}}
\put(3901,-3586){\makebox(0,0)[b]{\smash{$a$}}}
\put(3451,-1036){\makebox(0,0)[rb]{\smash{$z_1$}}}
\put(3451,-1936){\makebox(0,0)[rb]{\smash{$z_j$}}}
\put(3451,-2836){\makebox(0,0)[rb]{\smash{$z_N$}}}
\end{picture}
\caption{$B_{a+1,a-1}(t_1)\cdots B_{a+M, a-M}(t_M) \Omega_{a-M}(\vec z)$}
\label{fig:Psi}
\end{figure}

Graphically, a summand in \eqref{def:Psi} is denoted as in
\figref{fig:Psi}. 

\begin{rem}
The vector $\Psi(\vec t \,)$ gives an eigenvector of the transfer matrix
of the XYZ type spin chains when $\{ \nu, t_1, \ldots, t_M \}$ satisfies
so-called Bethe equations.  (See \cite{bax}, \cite{takh-fad}, \cite{take})
\end{rem}

The proof of the theorem is essentially the same as that of Theorem
1.4 of \cite{resh:92}. We first reduce the equation \eqref{diff-eq}
including the operator $A_j$ (\eqref{def:a-j:R-check} or
\eqref{def:a-j:R}) to the system
\begin{multline}
    \check R^{l_1, l_j}_{12}(z_1 - z_j - \kappa) 
    \ldots
    \check R^{l_{j-1}, l_j}_{j-1,j}(z_{j-1} - z_j - \kappa)
    f(\vec z_j)
    = \\
    =
    Z^{l_1, \ldots, l_{j-1}, l_{j+1}, \ldots, l_N, l_j}_{%
       z_1, \ldots, z_{j-1}, z_{j+1}, \ldots, z_N, z_j} (\kappa, c)
    \check R^{l_j, l_N}_{N-1,N}(z_j - z_N)
    \ldots
    \check R^{l_j, l_{j+1}}_{j,j+1}(z_j - z_{j+1})
    f(\vec z),
\label{diff-eq-2}
\end{multline}
using the unitarity \eqref{unitarity}.

\begin{figure}[ht]
\setlength{\unitlength}{0.00083300in}%
\begingroup\makeatletter\ifx\SetFigFont\undefined
% extract first six characters in \fmtname
\def\x#1#2#3#4#5#6#7\relax{\def\x{#1#2#3#4#5#6}}%
\expandafter\x\fmtname xxxxxx\relax \def\y{splain}%
\ifx\x\y   % LaTeX or SliTeX?
\gdef\SetFigFont#1#2#3{%
  \ifnum #1<17\tiny\else \ifnum #1<20\small\else
  \ifnum #1<24\normalsize\else \ifnum #1<29\large\else
  \ifnum #1<34\Large\else \ifnum #1<41\LARGE\else
     \huge\fi\fi\fi\fi\fi\fi
  \csname #3\endcsname}%
\else
\gdef\SetFigFont#1#2#3{\begingroup
  \count@#1\relax \ifnum 25<\count@\count@25\fi
  \def\x{\endgroup\@setsize\SetFigFont{#2pt}}%
  \expandafter\x
    \csname \romannumeral\the\count@ pt\expandafter\endcsname
    \csname @\romannumeral\the\count@ pt\endcsname
  \csname #3\endcsname}%
\fi
\fi\endgroup
\begin{picture}(4212,7224)(2701,-7273)
\thicklines
\put(5701,-361){\vector( 0,-1){2700}}
\multiput(5701,-3361)(0.00000,120.00000){3}{\line( 0, 1){ 60.000}}
\put(7501,-661){\vector(-1, 0){3000}}
\put(7501,-2461){\vector(-1, 0){3000}}
\put(7501,-1861){\line(-1, 0){2700}}
\put(4801,-1861){\vector( 0, 1){1500}}
\put(7501,-1261){\vector(-1, 0){3000}}
\multiput(5701,-4261)(0.00000,120.00000){3}{\line( 0, 1){ 60.000}}
\put(5701,-4261){\vector( 0,-1){2700}}
\put(7501,-4561){\vector(-1, 0){3000}}
\multiput(5101,-4261)(48.00000,0.00000){51}{\makebox(6.6667,10.0000){\SetFigFont{10}{12}{rm}.}}
\multiput(7501,-4261)(0.00000,-48.21429){57}{\makebox(6.6667,10.0000){\SetFigFont{10}{12}{rm}.}}
\multiput(7501,-6961)(-48.00000,0.00000){51}{\makebox(6.6667,10.0000){\SetFigFont{10}{12}{rm}.}}
\put(7501,-5761){\vector(-1, 0){3000}}
\put(7501,-5161){\vector(-1, 0){3000}}
\put(6901,-361){\vector( 0,-1){2700}}
\multiput(6901,-3361)(0.00000,120.00000){3}{\line( 0, 1){ 60.000}}
\put(6901,-4261){\vector( 0,-1){2700}}
\multiput(6901,-4261)(0.00000,120.00000){3}{\line( 0, 1){ 60.000}}
\multiput(7801,-661)(-120.00000,0.00000){3}{\line(-1, 0){ 60.000}}
\multiput(7801,-2461)(-120.00000,0.00000){3}{\line(-1, 0){ 60.000}}
\multiput(7801,-1861)(-120.00000,0.00000){3}{\line(-1, 0){ 60.000}}
\multiput(7801,-1261)(-120.00000,0.00000){3}{\line(-1, 0){ 60.000}}
\multiput(8101,-5761)(-109.09091,0.00000){6}{\line(-1, 0){ 54.545}}
\multiput(8101,-6061)(-120.00000,0.00000){3}{\line(-1, 0){ 60.000}}
\multiput(7801,-6061)(0.00000,120.00000){13}{\line( 0, 1){ 60.000}}
\multiput(7801,-4561)(-120.00000,0.00000){3}{\line(-1, 0){ 60.000}}
\multiput(8101,-5161)(-109.09091,0.00000){6}{\line(-1, 0){ 54.545}}
\put(7501,-6361){\vector(-1, 0){3000}}
\multiput(8101,-6361)(-109.09091,0.00000){6}{\line(-1, 0){ 54.545}}
\multiput(5701,-7261)(0.00000,120.00000){3}{\line( 0, 1){ 60.000}}
\multiput(6901,-7261)(0.00000,120.00000){3}{\line( 0, 1){ 60.000}}
\multiput(5101,-361)(48.00000,0.00000){51}{\makebox(6.6667,10.0000){\SetFigFont{10}{12}{rm}.}}
\multiput(7501,-361)(0.00000,-48.21429){57}{\makebox(6.6667,10.0000){\SetFigFont{10}{12}{rm}.}}
\multiput(7501,-3061)(-48.00000,0.00000){51}{\makebox(6.6667,10.0000){\SetFigFont{10}{12}{rm}.}}
\multiput(5701,-361)(0.00000,120.00000){3}{\line( 0, 1){ 60.000}}
\multiput(6901,-361)(0.00000,120.00000){3}{\line( 0, 1){ 60.000}}
\put(5776,-1036){\makebox(0,0)[lb]{\smash{$t_1$}}}
\put(5401,-3286){\makebox(0,0)[b]{\smash{$a$}}}
\put(5401,-4111){\makebox(0,0)[b]{\smash{$a$}}}
\put(5101,-4786){\makebox(0,0)[b]{\smash{$z_j+\kappa$}}}
\put(5101,-5386){\makebox(0,0)[b]{\smash{$z_1$}}}
\put(5101,-5986){\makebox(0,0)[b]{\smash{$z_{j-1}$}}}
\put(5776,-4936){\makebox(0,0)[lb]{\smash{$t_1$}}}
\put(6301,-1636){\makebox(0,0)[b]{\smash{$\cdots$}}}
\put(6301,-3286){\makebox(0,0)[b]{\smash{$\cdots$}}}
\put(6301,-2686){\makebox(0,0)[b]{\smash{$\cdots$}}}
\put(6301,-2236){\makebox(0,0)[b]{\smash{$\cdots$}}}
\put(6301,-1036){\makebox(0,0)[b]{\smash{$\cdots$}}}
\put(6301,-4111){\makebox(0,0)[b]{\smash{$\cdots$}}}
\put(6301,-5536){\makebox(0,0)[b]{\smash{$\cdots$}}}
\put(6301,-4936){\makebox(0,0)[b]{\smash{$\cdots$}}}
\put(6976,-1036){\makebox(0,0)[lb]{\smash{$t_M$}}}
\put(6976,-3286){\makebox(0,0)[lb]{\smash{$a+M$}}}
\put(6976,-4111){\makebox(0,0)[lb]{\smash{$a-M$}}}
\put(6976,-4936){\makebox(0,0)[lb]{\smash{$t_M$}}}
\put(5251,-1486){\makebox(0,0)[b]{\smash{$z_{j-1}$}}}
\put(5251,-886){\makebox(0,0)[b]{\smash{$z_1$}}}
\put(5251,-2086){\makebox(0,0)[b]{\smash{$z_j+\kappa$}}}
\put(5251,-2686){\makebox(0,0)[b]{\smash{$z_N$}}}
\put(6301,-6136){\makebox(0,0)[b]{\smash{$\cdots$}}}
\put(5101,-6586){\makebox(0,0)[b]{\smash{$z_N$}}}
\put(6301,-6736){\makebox(0,0)[b]{\smash{$\cdots$}}}
\put(6301,-7186){\makebox(0,0)[b]{\smash{$\cdots$}}}
\put(6976,-7186){\makebox(0,0)[lb]{\smash{$a+M$}}}
\put(5401,-7186){\makebox(0,0)[b]{\smash{$a$}}}
\put(5401,-211){\makebox(0,0)[b]{\smash{$a$}}}
\put(6301,-211){\makebox(0,0)[b]{\smash{$\cdots$}}}
\put(6976,-211){\makebox(0,0)[lb]{\smash{$a-M$}}}
\put(3901,-5611){\makebox(0,0)[lb]{\smash{$=$}}}
\end{picture}
\caption{Equivalence of \eqref{diff-eq-2} and \eqref{diff-eq-3}.}
\label{fig:prf-dif-eq3}
\end{figure}
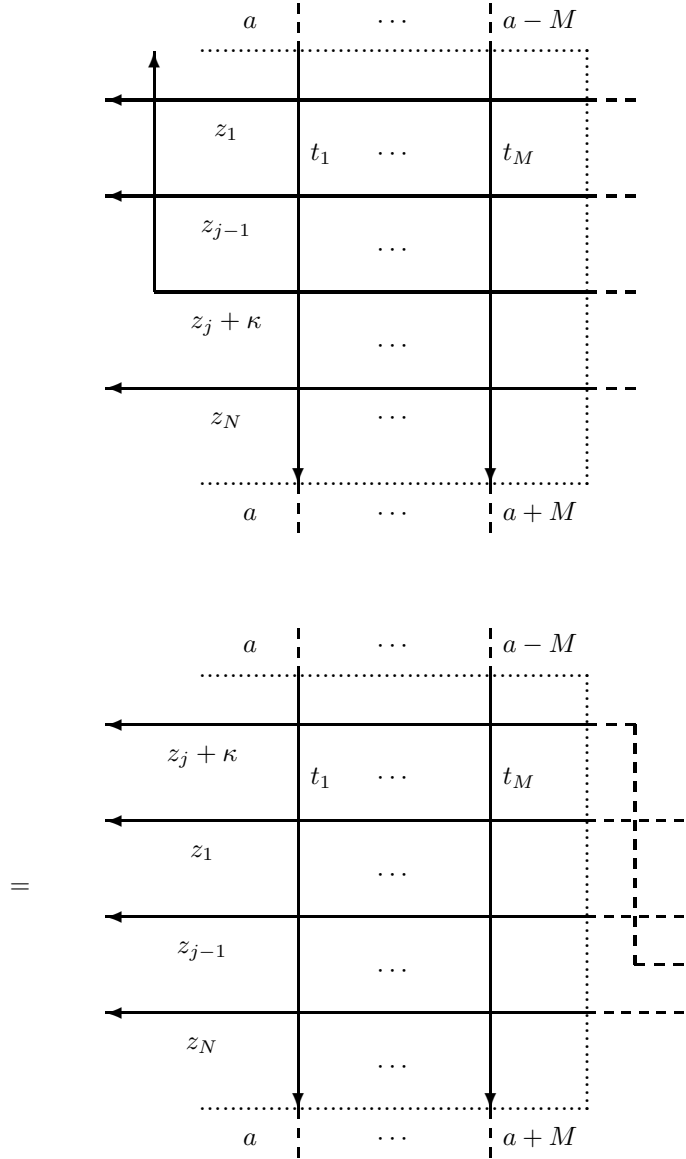

Assuming the form of the solution \eqref{solution}, we can show
that \eqref{diff-eq-2} is equivalent to
\begin{equation}
\begin{split}
    \sum_{a, \vec t} \varphi(\vec z_j | \vec t \,) &
    e^{2\pi i a \nu \eta}
    B^{l_N, \ldots, l_1, l_j}_{a+1,a-1}(t_1;z_N, \ldots, z_1, z_j + \kappa)
    \cdots \\
    \cdots &
    B^{l_N, \ldots, l_1, l_j}_{a+M,a-M}(t_M;z_N, \ldots, z_1, z_j + \kappa)
    \Omega^{l_j, l_1, \ldots, l_N}_{a-M}(z_j + \kappa, z_1, \ldots, z_N)
    = \\
    =
    \sum_{a, \vec t} \varphi(\vec z | \vec t \,) &
    Z\bigl(
    e^{2\pi i a \nu \eta}
    B^{l_j, l_N, \ldots, l_1}_{a+1,a-1}(t_1;z_j, z_N, \ldots, z_1)
    \cdots \\
    \cdots &
    B^{l_j, l_N, \ldots, l_1}_{a+M,a-M}(t_M;z_j, z_N, \ldots, z_1)
    \Omega^{l_1, \ldots, l_N, l_j}_{a-M}(z_1, \ldots, z_N, z_j)
    \bigr).
\end{split}
\label{diff-eq-3}
\end{equation}
Here, for example, $(z_N, \ldots, z_1, z_j)$ means
$(z_N, \ldots, z_{j+1}, z_{j-1}, \ldots, z_1, z_j)$.
Equivalence of the left hand side of \eqref{diff-eq-2} and
\eqref{diff-eq-3} is proved as illustrated in \figref{fig:prf-dif-eq3}.
We move the line with the spectral parameter $z_j + \kappa$, repeatedly
using the procedures \figref{fig:yb}, and then \figref{fig:R-on-int}. 
Because of the admissibility condition \eqref{admissible} and 
the equation \eqref{R-on-vac}, the IRF-type weights appearing in
\figref{fig:prf-dif-eq3} (crossings of dashed lines) are equal to the
unity. The right hand side of \eqref{diff-eq-3} is derived from the right
hand side of \eqref{diff-eq-2} in the same manner.

\newpage
In the next step we need a counterpart of Lemma 2.3 of \cite{resh:92}.
\begin{lem}
$(\text{Two-side formula})$ 

Fix $n$ $(1 \leqq n \leqq N-1)$ and $m$ $(1 \leqq m \leqq M)$ and let $a$
be an integer, $a_0 = a$ and $a_j = a + 2(l_1 + \cdots + l_j) \eta$ for
all $j \in \{1, \ldots, N\}$.  Then
\begin{equation}
\begin{split}
    &
    B^{l_N, \ldots, l_1}_{a_N-m+1, a_0+m-1}(t_m; z_N, \ldots, z_1) \cdots
    B^{l_N, \ldots, l_1}_{a_N, a_0}(t_1; z_N, \ldots, z_1)
    \Omega^{l_1, \ldots, l_N}_a(\vec z) =
\\
    = &
    \sum_{\{1, \ldots, m\} = \I \sqcup \II}
    \prod_{i \in \I, i' \in \II} \alpha(t_i, t_{i'})
    \prod_{i  \in \I}  \prod_{k = n+1}^N \alpha^{l_k}(t_i    - z_k)
    \prod_{i' \in \II} \prod_{k = 1}^{n} \delta^{l_k}(t_{i'} - z_k)
    \\
    &
    B^{l_n, \ldots, l_1}_{a_n - \s(\I) + \s(\II) + 1, a_0 + m - 1 }
     (t_{i_1}       ; z_n, \ldots, z_1)
    \cdots
    \\
    & \cdots
    B^{l_n, \ldots, l_1}_{a_n          + \s(\II)    , a_0 + \s(\II)}
     (t_{i_{\s(\I)}}; z_n, \ldots, z_1)
    \Omega^{l_1, \ldots, l_n}_{a_0 + \s(\II)}(z_1, \ldots, z_n)
    \tensor
    \\
    \tensor &
    B^{l_N, \ldots, l_{n+1}}_{a_N - m + 1 , a_n - \s(\I) + \s(\II) - 1 }
     (t_{i'_1}        ; z_N, \ldots, z_{n+1})
    \cdots
    \\
    & \cdots
    B^{l_N, \ldots, l_{n+1}}_{a_N - \s(\I), a_n - \s(\I)}
     (t_{i'_{\s(\II)}}; z_N, \ldots, z_{n+1})
    \Omega^{l_{n+1}, \ldots, l_N}_{a_n - \s(\I)}(z_{n+1}, \ldots, z_N),
\end{split}
\label{two-side-formula}
\end{equation}
where $\{1, \ldots, m\} = \I \sqcup \II$ is a partition, $\s$ denotes the
number of elements, $i_\bullet$ designates an element of $\I$ and
$i'_\bullet$ an element of $\II$.  

Due to \eqref{comm-rel-B-B}, the right hand side of
\eqref{two-side-formula} is well-defined.
\end{lem}

One can prove this lemma by induction on $m$. Using
\eqref{comm-rel-A-B}, \eqref{comm-rel-D-B} and a formula
\begin{equation*}
\begin{split}
    B^{l_N, \ldots, l_1    }_{a,a'}&(u;z_N ,\ldots, z_1)
    = \\ =&
    A^{l_N, \ldots, l_{n+1}}_{a,b }(u;z_N, \ldots, z_{n+1})
    B^{l_n, \ldots, l_1    }_{b,a'}(u;z_n, \ldots, z_1)
    + \\ +&
    B^{l_N, \ldots, l_{n+1}}_{a,b }(u;z_N, \ldots, z_{n+1})
    D^{l_n, \ldots, l_1    }_{b,a'}(u;z_n, \ldots, z_1)
\end{split}
\end{equation*}
for any integer $b$, which is easily derived from
\eqref{def:twisted-monodromy}, one can apply the formula (5.6) of
\cite{takh-fad} to this case.

Applying \eqref{two-side-formula} for $(m,n) = (M, 1)$ and $(M, N-1)$ to
the equations \eqref{diff-eq-3}, using \eqref{boundary-on-bethe} and
comparing the coefficients of
\begin{multline*}
    e^{2\pi i (a + 2 \s(\II) - 2l_j) \nu}
    B^{l_j}_{a + 1      , a + 2 \s(\II) - 2\l_j - 1}
      (t_{i'_1} + \kappa ;z_j + \kappa)
    \cdots
    \\
    \cdots
    B^{l_j}_{a + \s(\II), a +   \s(\II) - 2\l_j }
      (t_{i'_{\s(\II)}} + \kappa ;z_j + \kappa)
    \Omega^{l_j}_{a + \s(\II)} (z_j + \kappa)
    \tensor
    \\
    \tensor
    B^{l_N, \ldots, l_{j+1}, l_{j-1}, \ldots, l_1
       }_{a + 2 \s(\II) - 2 l_j + 1, a - 1}
     (t_{i_1}; z_N \ldots, z_{j+1}, z_{j-1}, \ldots, z_1)
    \cdots
    \\
    \cdots
    B^{l_N, \ldots, l_{j+1}, l_{j-1}, \ldots, l_1
       }_{a + M - 2 l_j + \s(\II), a - \s(\I)}
     (t_{i_{\s(\I)}}; z_N \ldots, z_{j+1}, z_{j-1}, \ldots, z_1)
    \\
    \Omega^{l_1, \ldots, l_{j-1}, l_{j+1}, \ldots, l_N
       }_{a + M - 2 l_j + \s(\II)}
      (z_1 \ldots, z_{j-1}, z_{j+1}, \ldots, z_N)
\end{multline*}
of the resulting equation for a partition
$\{1, \ldots, M\} = \I \sqcup \II$,
we can show that \eqref{solution} gives a solution to
\eqref{diff-eq} if $\varphi(\vec z | \vec t \,)$ satisfies
\begin{multline}
    \varphi(\vec z_j | t_i + \kappa \delta_{i, \II})
    \prod_{i \in \I, i' \in \II} \alpha(t_i' + \kappa, t_i)
    \prod_{i' \in \II} \prod^N \begin{Sb}k = 1 \\ k \neq j \end{Sb}
         \alpha^{l_k}(t_{i'}+\kappa - x_k)
    \prod_{i  \in \I } \delta^{l_j}(t_i - x_j - \kappa)
    =
\\
    =
    \varphi(\vec z | \vec t \,)
    e^{-2c \s(\II)}
    \prod_{i \in \I, i' \in \II} \alpha(t_i, t_i')
    \prod_{i  \in \I } \alpha^{l_j}(t_i - x_j)
    \prod_{i' \in \II} \prod^N \begin{Sb}k = 1 \\ k \neq j \end{Sb}
         \delta^{l_k}(t_{i'} - x_k).
\label{diff-eq-varphi}
\end{multline}
Here in the left hand side $\delta_{i,\II} = 0$ if $i \in \I$ and $=1$ if
$i \in \II$. The function defined by \eqref{def:varphi} is a solution of
\eqref{diff-eq-varphi}. Thus we have proved the theorem.

\subsection*{Ordinary difference equations and cycles}

We now return to the equations \eqref{diff-eq:F-l} and
\eqref{diff-eq:Phi} and give the solution to them in the form of infinite
products. 

First let us consider equation \eqref{diff-eq:F-l}.
Using the infinite product formula of the theta functions (see, e.g.,
\cite{mum}), we have
\begin{equation}
    \frac{\delta^l(t)}{\alpha^l(t + \kappa)}
    = e^{-4\pi i l \eta}
    \frac{(         e^{-2 \pi i t + 4 \pi i l \eta};p)_\infty}
         {(  q^{-1} e^{-2 \pi i t - 4 \pi i l \eta};p)_\infty}
    \frac{(p        e^{ 2 \pi i t - 4 \pi i l \eta};p)_\infty}
         {(p q      e^{ 2 \pi i t + 4 \pi i l \eta};p)_\infty},
\label{inf-prod:delta/alpha}
\end{equation}
where $p = \exp(2\pi i \tau)$, $q = \exp(2\pi i \kappa)$ and 
$(x;p)_\infty = \prod_{m=0}^\infty (1 - p^m x)$.
It is easy to see from \eqref{inf-prod:delta/alpha} that
\begin{equation}
    F^l(t) =
    e^{-4 \pi i l \eta t/\kappa}
    \frac{( q e^{-2\pi it + 4\pi il\eta}; p, q)_\infty}
         {(   e^{-2\pi it - 4\pi il\eta}; p, q)_\infty}
    \frac{(pq e^{ 2\pi it + 4\pi il\eta}; p, q)_\infty}
         {(p  e^{ 2\pi it - 4\pi il\eta}; p, q)_\infty}
\label{def:F-l}
\end{equation}
satisfies \eqref{diff-eq:F-l}, at least formally, where 
$(x;p,q)_\infty = \prod_{m=0}^\infty \prod_{n=0}^\infty (1 - p^m q^n x)$.

Similarly, equation \eqref{diff-eq:Phi} has a formal solution
\begin{equation}
    \Phi(t) =
    e^{4 \pi i \eta t/\kappa}
    (e^{-2\pi it}; p)_\infty (p e^{2\pi it}; p)_\infty 
    \frac{( q e^{-2\pi it - 4\pi i\eta}; p, q)_\infty}
         {(   e^{-2\pi it + 4\pi i\eta}; p, q)_\infty}
    \frac{(pq e^{ 2\pi it - 4\pi i\eta}; p, q)_\infty}
         {(p  e^{ 2\pi it + 4\pi i\eta}; p, q)_\infty}.
\label{def:Phi}
\end{equation}

In fact, the infinite products in \eqref{def:F-l} and \eqref{def:Phi}
gives meromorphic functions on the whole complex plane, provided that 
$\Im \kappa > 0$. The function $F^l(t)$ has zeros at
\begin{multline}
    \{t = n + m_1 \tau + m_2 \kappa + 2l\eta \, |\,
    n \in \Integer, m_1 \in \Integer_{\geq 0}, m_2 \in \Integer_{>0} \}
    \cup \\
    \cup
    \{t = n - m_1 \tau - m_2 \kappa - 2l\eta \, |\,
    n \in \Integer, m_1 \in \Integer_{> 0}, m_2 \in \Integer_{>0} \},
\label{zeros:F-l}
\end{multline}
and poles at
\begin{multline}
    \{t = n + m_1 \tau + m_2 \kappa - 2l\eta \, |\,
    n \in \Integer, m_1 \in \Integer_{\geq 0}, m_2 \in \Integer_{\geq 0} \}
    \cup \\
    \cup
    \{t = n - m_1 \tau - m_2 \kappa + 2l\eta \, |\,
    n \in \Integer, m_1 \in \Integer_{> 0}, m_2 \in \Integer_{\geq 0} \}.
\label{poles:F-l}
\end{multline}

Taking these properties of $F^l(t)$ into account, we can choose an
$N$-cycle in \eqref{solution} so that the sum over $t \in \calC$ reduces
to a finite sum, if $\kappa$ satisfies a certain rationality
condition. Suppose that there exist integers 
$(n, m_0, m_1) \in \Integer{>0} \times \Integer \times \Integer{>0}$
satisfying
\begin{equation}
    n\kappa = m_0 + m_1 \tau + 4l_j\eta
\label{integrality}
\end{equation}
for a certain $j \in \{1, \ldots, N\}$.
Then all but finite points in the set
\begin{equation}
    \calC':=
    \{t = m_0^0 + m_1^0 \tau + m \kappa + 2 l_j \eta \,|\,
    m \in \Integer \}
\label{cycle}
\end{equation}
fall into the set of zeros of $F^l(t)$ \eqref{zeros:F-l} for any integers
$m_0^0$ and $m_1^0 \geqq 0$. We can choose a cycle $\calC$ so that
$t_n - x_j \in \calC'$ for all $n = 1, \ldots, M$ and, thanks to the zeros
of function $F^{l_j}$, the sum over $t \in \calC$ in \eqref{solution} is
essentially a finite sum. Hence equation \eqref{solution} gives an
analytic solution of the difference equation \eqref{diff-eq}.

\section{Comments}
\setcounter{equation}{0}

In recent years several difference equations with elliptic coefficients
related to the $q$-Knizhnik-Zamolodchikov equations have been proposed:

\begin{itemize}
\item
Etingoff's equation:
Etingoff \cite{eti} showed that a (modified) trace of certain
intertwining operators of representations of quantum affine universal
enveloping algebras, $U_q(\hat{\frak g})$ satisfies a difference equation
with elliptic coefficients.

\item
Jimbo-Miwa-Nakayashiki's equation:
Jimbo, Miwa and Nakayashiki \cite{jmn} found a difference equation
which should be satisfied by correlation functions of eight vertex models.

\end{itemize}

Unfortunately we do not know the relation of our difference equation
\eqref{diff-eq} with any one of above equations. It might be possible that
one system turns into another by specialization of parameters. We can also
expect that a quasi-classical limit of our system \eqref{diff-eq} is
related to the Knizhnik-Zamolodchikov-Bernard equation.

It is a challenging problem to give a representation theoretical
interpretation to the solution \eqref{solution} like that of the integral
solution of the ($q$-)KZ equations \cite{mat}, \cite{fei-fre-resh}. 

%
%
%    Appendix A.
%
%
%%%%%%%%%%%%%%%%%%%%%%%%%%%%%%%%%%%%%%%%%%%%%%%%%%%%%%%%%%%%%%%%%%%%%%%%
\appendix
\renewcommand{\theequation}{\thesection.\arabic{equation}}
\renewcommand{\thethm}{\thesection.\arabic{thm}}
\section{Review of the Sklyanin algebra}
\setcounter{equation}{0}
\label{skl-alg}
In this appendix we recall several facts on the Sklyanin algebra and
its representations from \cite{skl:82} and \cite{skl:83}.
We use notations in \cite{mum} for theta functions:
\begin{equation}
    \theta_{ab}(z;\tau) = \sum_{n\in\Integer}
       \exp\left( 
             \pi i \left( \frac{a}{2} + n \right)^2 \tau
           +2\pi i \left( \frac{a}{2} + n \right)
                   \left( \frac{b}{2} + z \right)
           \right),
\label{theta:def}
\end{equation}
where $\tau$ is a complex number such that $\Im \tau > 0$.
The Pauli matrices are defined as usual:
\begin{equation}
    \sig^0 = \begin{pmatrix} 1 & 0  \\  0  & 1 \end{pmatrix},\quad
    \sig^1 = \begin{pmatrix} 0 & 1  \\  1  & 0 \end{pmatrix},\quad
    \sig^2 = \begin{pmatrix} 0 & -i \\  i  & 0 \end{pmatrix},\quad
    \sig^3 = \begin{pmatrix} 1 & 0  \\  0  &-1 \end{pmatrix}.
\end{equation}

The {\em Sklyanin algebra}, $U_{\tau,\eta}(sl(2))$ is generated
by four generators $S^0$, $S^1$, $S^2$, $S^3$, satisfying 
the following relations:
\begin{equation}
    R_{12}(u-v) L_{13}(u) L_{23}(v) =
    L_{23}(v) L_{13}(u) R_{12}(u-v).
\label{RLL}
\end{equation}
Here $u$, $v$ are complex parameters, the {\em $L$ operator},
$L(u)$, is defined by
\begin{equation}
    L(u) = \sum_{a=0}^3 W_a^L(u) \sig^a \tensor S^a,
\label{def:L}
\end{equation}
where
\begin{gather*}
    W_0^L(u) 
    = \frac{\theta_{11}(u;\tau)}
           {2\theta_{11}(2\eta;\tau)\theta_{11}(\eta;\tau)},\qquad
    W_1^L(u)
    = \frac{\theta_{10}(u;\tau)}
           {2\theta_{11}(2\eta;\tau)\theta_{10}(\eta;\tau)},
\\
    W_2^L(u)
    = \frac{\theta_{00}(u;\tau)}
           {2\theta_{11}(2\eta;\tau)\theta_{00}(\eta;\tau)},\qquad
    W_3^L(u)
    = \frac{\theta_{01}(u;\tau)}
           {2\theta_{11}(2\eta;\tau)\theta_{01}(\eta;\tau)},
\end{gather*}
$R(u)=R(u;\tau)$ is {\em Baxter's $R$ matrix} defined by
\begin{equation}
    R(u) = \sum_{a=0}^3 W_a^R(u) \sig^a \tensor \sig^a,\qquad
    W_a^R(u) := \theta_{11}(2\eta; \tau) W_a^L(u + \eta).
\label{def:R}
\end{equation}
and indices $\{0,1,2\}$ denote the spaces on which operators act
non-trivially: for example,
$$
    R_{12}(u) = 
    \sum_{a=0}^3
    W_a^R(u) 1 \tensor \sig^a \tensor \sig^a,\qquad
    L_{13}(u) =
    \sum_{a=0}^3
    W_a^L(u) \sig^a \tensor 1 \tensor S^a.
$$

The above relation \eqref{RLL} contains $u$ and $v$ as parameters,
but the commutation relations among $S^a$ ($a=0, \dots, 3$) do not
depend on them:
\begin{equation}
    [S^\alpha,S^0    ]_- =
    -i J_{\alpha,\beta} [S^\beta,S^\gamma]_+, \qquad
    [S^\alpha,S^\beta]_- =
                      i [S^0,    S^\gamma]_+,
\label{comm_rel}
\end{equation}
where $(\alpha, \beta, \gamma)$ stands for any cyclic permutation of
(1,2,3), $[A,B]_\pm = AB\pm BA$, and
$J_{\alpha,\beta}=(W_\alpha^2-W_\beta^2)/(W_\gamma^2-W_0^2)$ depend on
$\tau$ and $\eta$ but not on $u$.

The {\em spin $l$ representation} of the Sklyanin algebra,
$\rho^{l}: U_{\tau,\eta}(sl(2)) 
          \to \End_{\Complex}(\Theta^{4l}_{00})$
is defined as follows: The representation space $V^l$ is
\begin{equation}
    V^l = \Theta^{4l+}_{00} :=
    \{f(y) \, |\,
     f(y+1) = f(-y) = f(y), f(y+\tau)=\exp^{-4l\pi i(2y+\tau)}f(y) \}.
\end{equation}
It is easy to see that $\dim V^l = 2l+1$.
The generators of the algebra act on this space as difference operators:
\begin{equation}
    (\rho^l(S^a) f)(y) =
    \frac{s_a(y-l\eta)f(y+\eta)-s_a(-y-l\eta)f(y-\eta)}
         {\theta_{11}(2y;\tau)},
\end{equation}
where
\begin{alignat*}{2}
    s_0(y) &=  \theta_{11}(\eta;\tau) \theta_{11}(2y;\tau),\qquad&
    s_1(y) &=  \theta_{10}(\eta;\tau) \theta_{10}(2y;\tau),\\
    s_2(y) &= i\theta_{00}(\eta;\tau) \theta_{00}(2y;\tau),\qquad&
    s_3(y) &=  \theta_{01}(\eta;\tau) \theta_{01}(2y;\tau).
\end{alignat*}
These representations reduce to the usual spin $l$ representations of
$U(sl(2))$ for $J_{\alpha\beta} \to 0$ ($\eta \to 0$).
In particular, in the case $l= 1/2$, $S^a$  are expressed by
the Pauli matrices $\sigma^a$:
Let us identify $\Theta^{2+}_{00}$ and $\Complex^2$ by
\begin{equation}
\begin{aligned}
    \theta_{00}(2y;2\tau)-\theta_{10}(2y;2\tau) 
    &\longleftrightarrow
    \begin{pmatrix} 1 \\ 0 \end{pmatrix},
\\
    \theta_{00}(2y;2\tau)+\theta_{10}(2y;2\tau) 
    &\longleftrightarrow
    \begin{pmatrix} 0 \\ 1 \end{pmatrix}.
\end{aligned}
\label{rep:identify}
\end{equation}
Under this identification $S^a$ have matrix forms
\begin{equation}
    \rho^{1/2}(S^a) = \theta_{11}(2\eta;\tau) \sigma^a.
\label{rep:pauli}
\end{equation}

%
%
%
%     References.
%
%%%%%%%%%%%%%%%%%%%%%%%%%%%%%%%%%%%%%%%%%%%%%%%%%%%%%%%%%%%%

%
\end{document}